\documentstyle[12pt,amssymb,amscd]{amsart}
\theoremstyle{plain}
 \newtheorem{thm}{Theorem}[section]
 \newtheorem{lem}[thm]{Lemma}
 \newtheorem{prop}[thm]{Proposition}
 \newtheorem{cor}[thm]{Corollary}
\theoremstyle{definition}
 \newtheorem{defn}{Definition}[section]
\theoremstyle{remark}
 \newtheorem{rem}{Remark}[section]
\newcommand{\Ext}{\operatorname{Ext}}
\newcommand{\Hom}{\operatorname{Hom}}
\newcommand{\codim}{\operatorname{codim}}
\newcommand{\im}{\operatorname{im}}
\newcommand{\Gal}{\operatorname{Gal}}

\newcommand{\rk}{\operatorname{rk}}
\newcommand{\tr}{\operatorname{tr}}

\newcommand{\chr}{\operatorname{ch}}

\newcommand{\NS}{\operatorname{NS}}
\newcommand{\coker}{\operatorname{coker}}
\newcommand{\Pic}{\operatorname{Pic}}
\newcommand{\depth}{\operatorname{depth}}
\newcommand{\pdim}{\operatorname{proj-dim}}
\newcommand{\length}{\operatorname{length}}
\newcommand{\Hilb}{\operatorname{Hilb}}
\newcommand{\Alb}{\operatorname{Alb}}
\newcommand{\Tor}{\operatorname{Tor}}
\newcommand{ \Spec}{\operatorname{Spec}}
\newcommand{\Quot}{\operatorname{Quot}}
\newcommand{ \Def}{\operatorname{Def}}

\font\b=cmr10 scaled \magstep5
\def\bigzerou{\smash{\lower1.7ex\hbox{\b 0}}}

\numberwithin{equation}{section}
\begin{document}

\title{Some notes on the moduli of stable sheaves on elliptic surfaces}
\author{K\={o}ta Yoshioka\\
 Dept. of Math. Hiroshima Univ.}
 \address{Department of Mathematics, Faculty of Science, Hiroshima University,
Higashi-Hiroshima, 739, Japan }
\email{yoshioka@@top2.math.sci.hiroshima-u.ac.jp}
\thanks{Partially supported by JSPS Research Fellowships for Young Scientists}
 \maketitle

\section{Introduction}
Let $X$ be a smooth projective surface over $\Bbb C$
and $H$ an ample divisor on $X$.
Let $M_H(r,c_1,\Delta)$ be the moduli of stable sheaves $E$
of rank $r$ on $X$ with $c_1(E)=c_1 \in \NS(X)$ and $\Delta(E)=\Delta$,
where $\Delta(E):=c_2(E)-\{(\rk(E)-1)/2\rk(E)\} (c_1(E)^2)$.
In this note, we shall consider the moduli spaces on elliptic surfaces.
Let $\pi:X \to C$ be an elliptic surface such that every singular
fibre is irreducible and $f$ a fibre of $\pi$.
We assume that $X$ is regular, $(c_1,f)$ is odd
and $H$ is sufficiently close to $f$.
Then Friedman [F] showed that $M_H(2,c_1,\Delta)$ is birationally equivalent
to $S^n (J^d X)$,
where $n=\dim M_H(2,c_1,\Delta)/2$, $2d+1=(c_1,f)$ and 
$J^d X$ is an elliptic surface over $C$
whose generic fibre is the set of line bundles of degree $d$.
In this note, we shall generalize it to the case where 
$r$ and $(c_1,f)$ are relatively prime.

As an application, we shall show that 
$M_H(r,kH,\Delta)$ is a rational variety 
for the case where $(X,H) =(\Bbb P^2, \cal O_{\Bbb P^2}(1))$ and $(r,3k)=1$. 
We also consider moduli spaces on Abelian surfaces.
In particular, we shall compute a generator of 
$H^2(M_H(r,c_1,\Delta),\Bbb Z)$.
For general surfaces, Li [Li1, Li2] considered the structure of 
$H^i(M_H(2,c_1,\Delta),\Bbb Q)$, $i \leq 2$ and 
$\Pic(M_H(2,c_1,\Delta))\otimes \Bbb Q$
for $\Delta \gg 0$.
For the integral cohomologies, 
Mukai [Mu3, Mu5] and O'Grady [O] investigated the structure of 
$H^2(M_H(r,c_1,\Delta),\Bbb Z)$ and the Picard group, if $X$ is a K3 surface.
By the same method as in [Y2],
we get a generator of $H^2(M_H(r,c_1,\Delta),\Bbb Z)$,
if $X$ is a ruled surface.
Our results for Abelian surfaces are similar to these results.

In section 1, we shall consider the birational structure of 
$M_H(r,c_1,\Delta)$.
Our method is the same as that in Friedman [F] and Maruyama [M2].
That is, we shall use elementary transformations.
For simplicity, we assume that $X$ is regular.
Let $E$ be an element of $M_H(r,c_1,\Delta)$.
Since $H$ is sufficiently close to the fibre,
$E_{|\pi^{-1}(\eta)}$ is a stable vector bundle on $\pi^{-1}(\eta)$.
Then there is a stable vector bundle $E_1$ such that $E_{1|l}$ is semi-stable
in the sense of Simpson [S] for all fibres $l$,
and $E$ is obtained from $E_1$ by successive elementary transformations along
coherent sheaves of pure dimension 1 on fibres.
Let $E_2$ be a stable vector bundle such that 
$E_{2|\pi^{-1}(\eta)} \cong E_{1|\pi^{-1}(\eta)}$, 
$E_{2|l}$ is semi-stable
in the sense of Simpson and $\det E_{2|l} \cong \det E_{1|l}$ 
for all fibres $l$.
By using the irreducibility of $l$, we shall show that 
$E_2 \cong E_1 \otimes \pi^* L$, where $L \in \Pic(C)$.
Then we can easily show that $S^n (J^d X)$ is birationally equivalent to 
an irreducible component of $M_H(r,c_1,\Delta)$,
where $n=\dim M_H(r,c_1,\Delta)/2$ and $d$ is an integer.
By the dimension counting of non-locally free part 
(cf. [Y1, Thm. 0.4]), we see that every irreducible component contains
vector bundles (the non-locally free part is of codimension $r-1$).
Let $E$ be a vector bundle of $M_H(r,c_1,\Delta)$. 
We note that $\Ext^2(E,E(-l))_0 \cong \Hom(E,E(K_X+l))^{\vee}_0=0$
for all fibres $l$, where $\Ext^i(E,E(D))_0$ is the trace free part of 
$\Ext^i(E,E(D))$.
Then $\Ext^1(E,E)_0 \to \Ext^1(E_{|l},E_{|l})_0$ is surjective.
Considering the deformation space of $E_{|l}$, we shall show that
$S^n (J^d X)$ is birationally equivalent to $M_H(r,c_1,\Delta)$.

In section 2, we shall treat the moduli spaces on $\Bbb P^2$.
Let $V \subset H^0(\Bbb P^2,K_{\Bbb P^2}^{\vee})$ be a linear pencil
which contains an elliptic curve $C$.
Since $(K_{\Bbb P^2},H)<0$,
we can deform $E \in M_H(r,c_1,\Delta)$ to a sheaf $E' \in M_H(r,c_1,\Delta)$
such that $E'_{|C}$ is semi-stable.
If $(c_1,H)$ and $r$ are relatively prime,
then $E'_{|C}$ is a stable vector bundle.
Let $\Bbb P^2 \to \Bbb P^1$ be the rational map defined by $V$ and
$Y \to \Bbb P^2$ the blow-ups of $\Bbb P^2$ which defines the morphism
$Y \to \Bbb P^1$. 
Then $M_H(r,c_1,\Delta)$ is birationally equivalent to
a component of a moduli space $M_{H'}(r,c_1,\Delta)$,
where $H'$ is an ample divisor on $Y$ which is sufficiently
close to the fibre in $\NS(Y)$.
Since $M_{H'}(r,c_1,\Delta)$ is birationally equivalent to 
a symmetric product of $Y$,
we get that $M_H(r,c_1,\Delta)$ is rational.
We also prove that the moduli of simple torsion free sheaves on
Del Pezzo surfaces are irreducible.

In section 3, we shall consider the moduli spaces on an Abelian surface.
We assume that $c_1 \mod r\NS(X)$ is a primitive element of $\NS(X)/r\NS(X)$.
Mukai [Mu1] gave a complete description of $M_H(r,c_1,\Delta)$
in the case where $\dim M_H(r,c_1,\Delta)=2$.
Hence we assume that $\dim M_H(r,c_1,\Delta) \geq 4$.
By using a quasi-universal family [Mu3],
we shall construct a generator of 
$H^i(M_H(r,c_1,\Delta),\Bbb Z)$ for $i=1,2$,
where $H$ is a general polarization (Theorem \ref{thm:H2}).
Our method is the same as in G\"{o}ttsche and Huybrechts [G-H],
that is, we shall deform $X$ to a product of elliptic curves.
Then $M_H(r,c_1,0)$ is isomorphic to $X$ and
$M_H(r,c_1,\Delta)$ is birationally equivalent to 
$X \times Hilb_X^{r\Delta}$.
Since both spaces have trivial canonical bundles,
there are closed subsets $Z_1 \subset M_H(r,c_1,\Delta)$
and $Z_2 \subset X \times Hilb_X^{r\Delta}$ such that
$\codim(Z_1) \geq 2$, $\codim(Z_2) \geq 2$ and 
$M_H(r,c_1,\Delta) \setminus Z_1 \cong (X \times Hilb_X^{r\Delta}) 
\setminus Z_2$.
Hence we get an isomorphism
$H^i(M_H(r,c_1,\Delta),\Bbb Z) \cong H^i(X \times Hilb_X^{r\Delta},
\Bbb Z)$, $i=1,2$.
Constructing a family of stable sheaves parametrized by 
$X \times Hilb_X^{r\Delta} \setminus Z_2$ directly,
we shall construct a generator of $H^i(M_H(r,c_1,\Delta),\Bbb Z)$,
$i=1,2$.
By using deformation of $X$ and the result in [Y4],
we shall also show that
the Betti numbers of $M_H(2,c_1,\Delta)$ are the same as those of
$M_H(1,0,2\Delta)$ (Theorem \ref{thm:B}).
We next show that the morphism $M_H(r,c_1,\Delta) \to \Pic^0(X) \times X$
defined in [Y2, Sect. 5] is an Albanese map, 
if $\dim M_H(r,c_1,\Delta) \geq 4$.
Combining all together, we also describe the Picard group of 
$M_H(r,c_1,\Delta)$ (Theorem \ref{thm:Pic}).
   
I would like to thank Professors A. Ishii and M. Maruyama
 for valuable discussions.

\vspace{ 1pc}

Notation.\newline
Let $X$ be a smooth projective surface over $\Bbb C$
and $H$ an ample divisor on $X$.
For a scheme $S$, we denote the projection $S \times X \to S$
by $p_S$.
We denote the N\'{e}ron-Severi group of $X$ by $\NS(X)$.
For an $x \in \NS(X) \otimes \Bbb Q$, we set 
$P(x):=(x,x-K_X)/2+\chi(\cal O_X)$.
  
For a torsion free sheaf $E$ on $X$, we set
$$
\Delta(E):=c_2(E)-\frac{\rk(E)-1}{2\rk (E)}(c_1(E)^2).
$$
We denote the trace free part of $\Ext^i(E,E(D))$ by $\Ext^i(E,E(D))_0$.

In this note, we only use the notion of
(semi-)stability in the sense of Mumford.
Let $M_H(r,c_1,\Delta)$ be the moduli of stable sheaves $E$
of rank $r$ on $X$ with $c_1(E)=c_1 \in \NS(X)$ and $\Delta(E)=\Delta$.
We denote the open subscheme of $M_H(r,c_1,\Delta)$ consisting
of stable vector bundles by $M_H(r,c_1,\Delta)_0$.

\section{Moduli spaces on elliptic surfaces}
\subsection{Preliminaries}
Let $\pi:X \to C$ be an elliptic surface
such that every fibre is irreducible.
We denote the algebraically equivalence class of a fibre by $f$.
Let $\eta$ be the generic point of the base curve $C$.
Let $J^d X \to C$ be the elliptic surface over $C$
such that the generic fibre is the set of line bundles of degree 
$d$ on $X_{|\pi^{-1}(\eta)}$.
For a coherent sheaf $E$ on a fibre $l$,
we set
\begin{align*}
\rk(E) &:=\length_{\cal O_{\eta_l}}(E \otimes \cal O_{\eta_l}),\\
\deg(E) &:=\chi(E),
\end{align*}
where $\eta_l$ is the generic point of $l$.

A coherent sheaf $E$ of pure dimension 1 on a fibre $l$ is semi-stable
if 
$$
\frac{\chi(F)}{\rk(F)} \leq \frac{\chi(E)}{\rk(E)}
$$
for all subsheaves $F \ne 0$ of $E$.  

\begin{lem}\label{lem:k}
Let $L$ be a relatively ample divisor on $X$.
Let $D$ be a divisor on $X$ such that 
$(D,f)\ne 0$ and $(D,L+kf)=0$ for some positive number $k$.
Then,
\begin{equation}\label{eq:2}
(D^2) \leq \frac{-1}{(L,f)^2}((L^2)+2k(L,f)).
\end{equation}
\end{lem}

\begin{pf}
We set $D=aL+bf+D'$, where $a, b \in \Bbb Q$
and $(D',L)=(D',f)=0$.
By the Hodge index theorem, $({D'}^2)\leq 0$.
Hence $(D^2)=((aL+bf)^2)+({D'}^2) \leq
((aL+bf)^2)=a^2(L^2)+2ab(L,f)$.
Thus we may assume that $D=aL+bf$.
$(D,L+kf)=0$ implies that $b(L,f)=-a(L,L+kf)$.
Hence $((aL+bf)^2)=-a^2((L^2)+2k(L,f))$.
Since $(L,f) \ne 0$, we get that $|a| \geq 1/|(L,f)|$.
Hence \eqref{eq:2} holds.
\end{pf}

\begin{lem}
Let $r$ be a positive integer and
$c_1$ an algebraically equivalence class on $X$ such that
$(c_1,f)$ and $r$ are relatively prime.
Let $L$ be an ample divisor on $X$.
Then 
\begin{equation*}
M_{L+nf}(r,c_1,\Delta)
=\left\{E\left|\begin{split} 
&\text{$E$ is torsion free of rank $r$
with $(c_1(E),\Delta(E))$}\\
&\text{$=(c_1,\Delta)$ and 
$E_{|\pi^{-1}(\eta)}$ is stable}
\end{split}
\right.
\right\}
\end{equation*}
for $n>(r^3(L,f)^2\Delta-2(L^2))/4(L,f)^2$.
We denote this space by $M(r,c_1,\Delta)$.
\end{lem}

\begin{pf}
The proof is similar to that in [Y3, Prop. 6.2]
(in [Y3], we used slightly different definition of $\Delta$).
\end{pf}
Since $\Ext^2(E,E)_0 \cong \Hom(E,E)_0^{\vee}=0$, $E \in M(r,c_1,\Delta)$,
$M(r,c_1,\Delta)$ is smooth of dimension $2r\Delta-
(r^2-1)\chi(\cal O_X)+\dim \Pic^0(X)$.
For a stable sheaf $E \in M(r,c_1,\Delta)$, $\chi(E_{|f})=(c_1,f)$ and $\chi(E \otimes k_x)=r$ are relatively prime,
where $E$ is locally free at $x \in X$ and $k_x$ is the structure sheaf of $x$.
Hence there is a universal family
(cf. [M1, Thm. 6.11]). 
If we fix the rank $r$ and the equivalence class $c_1 \mod \pi^*H^1(C,\Bbb Z)$,
then we may denote $M(r,c_1,\Delta)$ by $M(\Delta)$.
In fact, $c_1 \mod r\pi^* H^2(C,\Bbb Z)$ is determined by $r\Delta$
and the isomorphic class of $M_H(r,c_1,\Delta)$ is determined by
$r$, $c_1 \mod r\pi^* H^2(C,\Bbb Z)$ and $\Delta$. 
\begin{lem}\label{lem:D}
Let $E$ be a vector bundle of rank $r$ on $X$ such that
$(c_1(E),f)=d$,  
and let $F$ be a coherent sheaf of pure dimension 1
 on a fibre $l$
with $\rk(F)=r_1$ and $\deg(F)=d_1$.
Let $E \to F$ be a surjective homomorphism
and $E'$ the kernel.
Then 
\begin{equation}
\Delta(E')=\Delta(E)+\frac{rd_1-r_1d}{r}.
\end{equation}
\end{lem}

\begin{pf}
For a coherent sheaf $G$ on $X$,
$\chi(G)=\rk (G) P(c_1(G)/\rk G)-\Delta (G)$.
Since $\chi (E)=\chi (E')+\chi(F)$,
\begin{align*}
\Delta(E')-\Delta(E) & = d_1-r(P(c_1(E)/\rk E)-P(c_1(E')/\rk E'))\\
& = d_1-\frac{r_1 d}{r}.
\end{align*}
\end{pf}
The following is a special case of Maruyama [M2, 3.8].
\begin{cor}\label{cor:1}
Let $E$ be a vector bundle on $X$ such that $E_{|\pi^{-1}(\eta)}$
is a semi-stable vector bundle.
Then there is a vector bundle $E'$ on $X$
such that $E'_{|l}$ is semi-stable for every fibre $l$ and 
$E$ is obtained from $E'$ by successive elementary transformations
along coherent sheaves of pure dimension 1 on fibres.
\end{cor}

\begin{pf}
We note that $\Delta(E) \geq 0$.
We shall prove our claim by induction on $\Delta(E)$.
We assume that there is a fibre $l$ such that $E_{|l}$ is not semi-stable.
Then there is a surjective homomorphism
$E_{|l} \to F$ such that 
$F$ is of pure dimension 1 and $\chi(E_{|l})/\rk (E_{|l})>\chi(F)/\rk F$.
We shall consider the following elementary transformation
along $F$:
\begin{equation*}
0 \to E_1 \to E \to F \to 0.
\end{equation*}
Since $\depth_{\cal O_x}F=1$, $x \in C$
and $X$ is smooth, we see that $\pdim_{\cal O_x} F=\dim X-\depth_{\cal O_x}F=1$.
Hence $E_1$ is also locally free. 
By Lemma \ref{lem:D},
we get that $\Delta(E_1)<\Delta(E)$.
Hence we obtain our corollary.
\end{pf}

\subsection{General element of $M(\Delta)$}
Let $E$ be a general element of $M(\Delta)$.
We shall consider the Harder-Narasimhan filtration of the restriction
$E_{|l}$ of $E$ to fibres $l$.
In particular, we shall show that $E_{|l}$ is semi-stable for all singular
fibres $l$.
\begin{lem}\label{lem:def}
Let $C$ be a projective curve and $\cal O_C(1)$ an ample divisor on $C$.
Let $L$ be a line bundle on $C$.
Let $Q$ be the subscheme of $Quot_{\cal O_C(-n)^{\oplus N}/C}$ 
parametrizing quotients $\cal O_C(-n)^{\oplus N} \to E$
such that 
$(\mathrm{i})$ $E$ is a locally free sheaf of rank r with $\det E = L$ and
$(\mathrm{ii})$ $H^1(C,E(n))=0$.
Then $Q$ is smooth and irreducible.
\end{lem} 
\begin{pf}
Let $\lambda:\cal O_C(-n)^{\oplus N} \to E$ be a quotient which belongs to 
$Q$. Then we see that $\Ext^1(\ker \lambda,E)=0$.  
Since $\Hom(\ker \lambda,E) \to \Ext^1(E,E) \overset{\tr}\to H^1(C,\cal O_C)$
is surjective, 
$Q$ is smooth.
For $l \geq n$,
there is an exact sequence
$ 0 \to \cal O_C^{\oplus (r-1)} \to E(l) \to L(rl) \to 0$.
We set $\Bbb P:=\Bbb P(\Ext^1(L(rl),\cal O_C^{\oplus (r-1)})^{\vee})$.
We shall consider the universal extension:
$$
0 \to \cal O_{\Bbb P \times C}^{\oplus (r-1)} \to
 \cal E \to L(rl) \otimes \cal O_{\Bbb P}(-1) \to 0.
$$
Let $\Bbb P'$ be the open subscheme of $\Bbb P$
of points $y$ such that $H^1(C,\cal E_y)=0$.
Then $p_{\Bbb P'*}(\cal E)$ is a locally free sheaf on $\Bbb P'$.
Let $\phi:\Bbb A \to \Bbb P'$ be the vector bundle associated 
to the locally free sheaf 
$\cal Hom(\cal O_{\Bbb P'}^{\oplus N},p_{\Bbb P'*}(\cal E))$.
Then there is a homomorphism $\Lambda:\cal O_{\Bbb A \times C}^{\oplus N}
 \to (\phi \times 1)^* \cal E$.
Let $\Bbb A'$ be the open subscheme of $\Bbb A$ such that 
$\Lambda$ is surjective.
Then there is a surjective morphism
$\Bbb A' \to Q$, and hence $Q$ is irreducible.
\end{pf}

\begin{prop}
Let $M(\Delta)^0$ be the open subscheme of $M(\Delta)$ of elements $E$
such that $E_{|l}$ is semi-stable for every singular fibre $l$.
Then $M(\Delta)^0$ is a dense subscheme of $M(\Delta)$.
\end{prop}
\begin{pf}
Let $E$ be an element of $M(\Delta)$.
Since $E_{|\pi^{-1}(\eta)}$ is stable,
we see that $\Ext^2(E,E(-l))_0 \cong \Hom(E,E(l+K_X))_0^{\vee}=0$.
Hence we get that
$\Ext^1(E,E)_0 \to \Ext^1(E_{|l},E_{|l})_0$ is surjective.
Let $m$ be the multiplicity of $l$ and set $l=ml'$.
By Corollary \ref{cor:1}, there is a vector bundle $E_1$ on $X$ such that
$E_{1|l}$ is semi-stable and $\det(E_{1|l})=\det(E_{|l})\otimes \cal O_l(kl')$.
Since $(r,m)=1$,
replacing $E_1$ by $E_1 \otimes \cal O_X(jl')$, we may assume that
$\det(E_{1|l})=\det(E_{|l})$. 
By using Lemma \ref{lem:def}, 
we see that $E_{|l}$ deforms to a semi-stable vector bundle.
Hence we see that $E$ deforms to a sheaf $E'$ such that
$E'_{|l}$ is semi-stable.
Thus $M(\Delta)^0$ is an open dense subscheme of $M(\Delta)$. 
\end{pf}

\begin{lem}
Let $l$ be a smooth fibre.
Let $h:=\{(r_1,d_1),(r_2,d_2),\dots ,(r_s,d_s)\}$
be a sequence of pairs of integers such that $r_i>0$, $1 \leq i \leq s$
and $d_1/r_1>d_2/r_2>\dots>d_s/r_s$.
Let $D_h$ be the subset of $M(r,c_1,c_2)$ of elements $E$ 
such that the Harder-Narasimhan filtration of $E_{|l}$
: $0 \subset F_1 \subset F_2 \subset \dots \subset F_s=E_{|l}$
satisfies that $\rk(F_i/F_{i-1})=r_i$
and $\deg(F_i/F_{i-1})=d_i$, $1 \leq i \leq s$.
Then $\codim(D_h)\geq \sum_{i<j}r_jd_i-r_id_j$.
In particular, if $\codim(D_h)=1$,
then $s=2$ and $r_2d_1-r_1d_2=1$.
\end{lem}  

\begin{pf}
Let $\Def(E_{|l})$ be the local deformation space
of $E_{|l}$ of fixed determinant and
$\Def(E_{|l})_h$ the subset of $\Def(E_{|l})$ of elements $G$
such that the Harder-Narasimhan filtration of $G$
: $0 \subset F_1 \subset F_2 \subset \dots \subset F_s=G$
satisfies that $\rk(F_i/F_{i-1})=r_i$
and $\deg(F_i/F_{i-1})=d_i$, $1 \leq i \leq s$.
We assume that $\Def(E_{|l})_h$ is not empty.
We note that 
$\Ext^1(E,E)_0 \to \Ext^1(E_{|l},E_{|l})_0$ is surjective.
It is known that $\codim(\Def(E_{|l})_h)=\sum_{i<j}r_jd_i-r_id_j$
(cf. [A-B, Thm. 7.14]).
Hence we get our lemma.
\end{pf}  
Let $(r_1,d_1)$ be the pair of integers such that
$0<r_1<r$ and $rd_1-r_1d=1$.
Let $M(\Delta)^1$ be the open subscheme of $M(\Delta)^0$ of 
elements $E$ such that
$E_{|l}$ is stable, or
the Harder-Narasimhan filtration of $E_{|l}$ is
$0 \subset F \subset E_{|l}$ for all fibres $l$,
where $F$ is a stable vector bundle of rank $r_1$ on $l$
 with $\deg(F)=d_1$.
Then $M(\Delta)^1$ is an open dense subscheme of $M(\Delta)^0$.

\subsection{Vector bundles on elliptic curves}
The following is due to Atiyah [A].
\begin{lem}\label{lem:bdle}
Let $C$ be a smooth elliptic curve.
Let $r$ be a positive integer and $d$ an integer such that 
$(r,d)=1$. Then,\newline
$(1)$ There is a stable vector bundle of rank $r$ and degree $d$.
\newline
$(2)$ Let $(r_1,d_1)$ be the pair of integers such that 
$r_1d-rd_1=1$ and $0<r_1<r$.
Let $E_1$ be a stable vector bundle of rank $r_1$ and degree $d_1$.
Then every stable vector bundle $E$ of rank $r$ and degree $d$
is defined by an exact sequence
\begin{equation}
0 \to E_1 \to E \to E_2 \to 0,
\end{equation}
where $E_2$ is a stable vector bundle of rank $r_2:=r-r_1$
and degree $d_2:=d-d_1$.\newline
$(3)$ Let $0 \subset F_1 \subset F_2 \subset \dots \subset F_s=E$ be the 
Harder-Narasimhan filtration of a vector bundle $E$.
Then $E \cong \oplus_{i=1}^s E_i$, where $E_i:=F_i/F_{i-1}$. 
\end{lem}

\begin{pf}
(1) We shall prove our claim by induction on $r$.
If $r=1$, then our claim obviously holds.
Let $(r_1,d_1)$ be the pair of integers such that 
$r_1d-rd_1=1$ and $0<r_1<r$.
We set $r_2:=r-r_1$ and $d_2:=d-d_1$.
By induction hypothesis, there are stable vector bundles
$E_i$ of rank $r_i$ and degree $d_i$, $i=1,2$.
Since $d_1/r_1<d_2/r_2$,
$\Hom(E_2,E_1)=0$.
By using the Riemann-Roch theorem,
we get that $\Ext^1(E_2,E_1) \cong \Bbb C$.
Let $0 \to E_1 \to E \to E_2 \to 0$ be a non-trivial
extension.
We shall show that $E$ is stable.
If $E$ is not stable, then there is a semi-stable subsheaf $G$
of $E$ such that $\deg G/\rk G>d/r$.
Since $G$ and $E_2$ are semi-stable and $G \to E\to E_2$ is not zero,
 $\deg G/\rk G \leq d_2/r_2$.
We assume that $\deg G/\rk G<d_2/r_2$.
Then we see that 
$1/rr_2=d_2/r_2-d/r>d_2/r_2-\deg G/\rk G \geq 1/r_2\rk G$,
which is a contradiction.
Hence $\deg G/\rk G=d_2/r_2$.
Then we get that $\rk G=r_2$ and $\deg G=d_2$.
Hence $G \cong E_2$, which is a contradiction.

(2) Let $E$ be a stable vector bundle of rank $r$
and degree $d$.
Then $\Ext^1(F_1,E) \cong \Hom(E,F_1)^{\vee}=0$.
By the Riemann-Roch theorem, there is a non-zero
homomorphism $\varphi:E_1 \to E$.
We shall show that $\varphi$ is injective and 
$\coker \varphi$ is stable.
Since $E_1$ and $E$ are stable, $d_1/r_1 \leq 
\deg \varphi(E_1)/\rk  \varphi(E_1)<d/r$.
In the same way as in the proof of (1),
we see that $\rk \varphi(E_1)=r_1$ and
$\deg \varphi(E_1)=d_1$.
Hence we get that $E_1 \cong \varphi(E_1)$.
We set $E_2:=\coker \varphi$.
We assume that there is a quotient $G$ of $E_2$
such that
$G$ is semi-stable and $d_2/r_2>\deg G/\rk G$.
Since $G$ is a quotient of $E$,
we get that $d/r<\deg G/\rk G$.
Hence we get that $d/r<\deg G/\rk G<d_2/r_2$.
Then $1/rr_2=d_2/r_2-d/r>d_2/r_2-\deg G/\rk G \geq
1/r_2\rk G$,
which is a contradiction.
Hence $E_2$ is a stable vector bundle.

(3) Since $\deg E_i/\rk E_i>\deg E_j/\rk E_j$,
$i<j$, the Serre duality implies that $\Ext^1(E_j,E_i)=0$, $i<j$.
 By the induction on $s$, we see that $E \cong \oplus_i E_i$.
\end{pf}

\begin{lem}\label{lem:elm}
Let $(r,d)$ (resp. $(r_1,d_1)$, $(r_2,d_2)$) be the pair in 
Lemma \ref{lem:bdle}.
Let $E$ be a vector bundle of rank $r$ on $C$
with degree $d$ and $E_2$ a stable vector bundle of rank $r_2$ on $C$
with degree $d_2$.\newline
(1) If $E$ is stable, then $\Hom(E,E_2) \cong \Bbb C$
and a non-zero homomorphism is surjective.\newline
(2) Let $F_1$ (resp. $F_2$) be a stable vector bundle of 
rank $r_1$ and degree $d_1$ (resp. rank $r_2$ and degree $d_2$).
We assume that $E \cong F_1 \oplus F_2$
and there is a surjective homomorphism $\varphi:E \to E_2$
such that $\ker \varphi$ is also stable.
Then $E_2 \cong F_2$ and $\Hom(E,E_2) \cong \Bbb C^{\oplus 2}$.
\end{lem} 

\begin{pf}
(1) Since $E$ is stable, $\Ext^1(E,E_2) \cong \Hom(E_2,E)^{\vee}=0$.
By the Riemann-Roch theorem, we see that
$\dim \Hom(E,E_2)=1$.
In the same way as in the proof of Lemma \ref{lem:bdle},
we see that a non-zero homomorphism
$E \to E_2$ is surjective.

(2) If $E_2 \not\cong F_2$, then $\ker \varphi \cong 
\ker(\varphi_{|F_1})\oplus F_2$.
Since $\varphi_{|F_1}:F_1 \to E_2$ is surjective, 
$\ker(\varphi_{|F_1}) \ne 0$. 
Hence $E_2 \cong F_2$.
By the Riemann-Roch theorem, $\Hom(F_1,E_2) \cong \Bbb C$.
Therefore $\Hom(E,E_2) \cong \Bbb C^{\oplus 2}$.
\end{pf}

Let $C_0$ be the open subscheme of $C$ such that 
$\pi:X_0:=X \times_C C_0 \to C_0$ is smooth.
We assume that $\pi$ has a section $\sigma$.
We denote the relative moduli space of stable vector bundles
of rank $r$ on fibres with degree $d$
by $\cal M_{X_0/C_0}(r,d) \to C_0$.
We assume that $(r,d)=1$.
We shall construct a family of stable vector bundles
$\cal E_{r,d}$ on $X_0 \times_{C_0}X_0$
and show that $\cal M_{X_0/C_0}(r,d) \cong X_0$ as a $C_0$-scheme,
by using induction on $r$.
If $r=1$, then
$\cal E_{1,d}:=\cal O_{X_0 \times _{C_0}X_0}((d+1)\sigma-\Delta)$
is a universal family, where
$\Delta$ is the diagonal of $X_0 \times _{C_0}X_0$.
Let $(r_1,d_1)$ be the pair of integers
such that $r_1d-rd_1=1$ and $0<r_1<r$.
We set $r_2=r-r_1$ and $d_2=d-d_1$.
Let $E$ be a vector bundle on $X_0$
such that $E_{|l}$ is a stable vector bundle
of rank $r_2$  and $\det E_{|l} \cong \cal O_l(d_2\sigma)$ 
for every fibre $l$. 
By using Lemma \ref{lem:bdle}, we see that
$\cal L:=\Ext^1_{p_{X_0}}(E,\cal E_{r_1,d_1})$
 is a line bundle on $X_0$.
Then there is the universal extension
\begin{equation}
0 \to \cal E_{r_1,d_1} \to \cal E_{r,d} \to 
E \otimes p_{X_0}^*(\cal L) \to 0,
\end{equation}
which parametrizes stable vector bundles of rank $r$
on fibres with degree $d$.
Hence there is a morphism
$X_0 \to \cal M_{X_0/C_0}(r,d)$.
By our construction, this morphism is injective.
By ZMT, it is an isomorphism.
 
\begin{lem}\label{lem:hom}
Let $E$ and $E'$ be semi-stable vector bundles 
on a multiple fibre $l=ml'$
such that $\rk E=\rk E'$,
$\det E \cong \det E'$, and $\chi(E)=\chi(E')=d$.
Then,
\begin{equation}
\Hom(E,E')=
\begin{cases}
\Bbb C,& \text{ if $E \cong E'$},\\
 0,& \text{ otherwise}.
\end{cases}
\end{equation}
\end{lem}

\begin{pf}
We set $L:=\cal O_X(-l')_{|l'}$.
We note that $\rk(E \otimes L^{\otimes k})=r$ and 
$\chi(E \otimes L^{\otimes k})=d/m$ for $0 \leq k \leq m-1$.
Since $(r,d)=1$ and $E$ is semi-stable,
$E \otimes L^{\otimes k}$ is a stable sheaf on $ml'$.
Thus 
$0 \subset E(-(m-1)l') \subset E(-(m-2)l') \subset
\dots \subset E(-l') \subset E$
is a Jordan-H\"{o}lder filtration of $E$.
Since the order of $L \in Pic^0(l')$ is $m$ and $(m,r)=1$,
$\det E \cong \det E'$ and the stabilities of 
$E_{|l'}$ and $E_{|l'}$ imply that 
$\Hom(E_{|l},E'\otimes L^{\otimes k})=0$ for
$1 \leq k \leq m-1$.
Let $\varphi:E \to E'$ be a non-zero homomorphism.
We shall show that $\varphi$ is an isomorphism.
Since $\Hom(E_{|l},E'\otimes L^{\otimes k})=0$ for
$1 \leq k \leq m-1$, 
we see that $\varphi_{|l'}:E_{|l'} \to E'_{|l'}$ is not zero,
which implies that $E_{|l'} \cong E'_{|l'}$.
By Nakayama's lemma, $\varphi$ is an isomorphism.
Then it is easy to see that 
$\Hom(E,E') \cong \Bbb C$.
\end{pf}

\begin{lem}\label{lem:5-10}
Let $E,E'$ be vector bundles of rank $r$ on $X$
such that $E_{|l}$ and $E'_{|l}$
are semi-stable for all fibres $l$ and 
$\det E \cong \det E'$.
Then there is a line bundle $L$ on $C$
such that $E \cong E' \otimes \pi^*(L)$.
\end{lem}
\begin{pf}
We note that $E_{|\pi^{-1}(\eta)} \cong E'_{|\pi^{-1}(\eta)}$.
By the upper semi-continuity of 
$h^0(l,E^{'\vee}\otimes E_{|l})$,
there is a non-zero homomorphism $E'_{|l} \to E_{|l}$
for every fibre $l$.
Since $E_{|l}$ and $E'_{|l}$
are semi-stable,
Lemma \ref{lem:hom} implies that
$E_{|l} \cong E'_{|l}$ and $H^0(l,E^{'\vee}\otimes E_{|l})
\cong \Bbb C$.
By the base change theorem, we get that
$L:=\pi_*(E^{'\vee} \otimes E)$ is a line on $C$
and $\pi^*(L) \otimes E' \to E$ is an isomorphism.
\end{pf}

\begin{cor}
$M(\Delta)$ is not empty
if and only if
$\Delta\geq \Delta_0:=\frac{(r^2-1)}{2r} \chi(\cal O_X)$.
\end{cor}

\begin{pf}
We set $\Delta':=\min\{\Delta|M(\Delta) \ne \emptyset \}$.
Lemma \ref{lem:5-10} implies that 
$\dim \Pic^0(X)=\dim M(\Delta')=2r\Delta'-(r^2-1)\chi(\cal O_X)
+\dim \Pic^0(X)$.
Hence we get our claim
\end{pf}

\begin{rem}
Let $E$ be an element of $M(\Delta_0)$.
By Lemma \ref{lem:5-10},
there is a surjective morphism
$\Pic^0(X) \to M(\Delta_0)$
sending $L \in \Pic^0(X)$ to $E \otimes L$.
Hence we get that $M(\Delta_0)=\Pic^0(X)/\Phi(E)$,
where $\Phi(E):=\{L \in \Pic^0(X)|E \otimes L \cong E \}$.
In particular, if $\Pic^0(X)=\Pic^0(C)$, then 
$M(\Delta_0)=\Pic^0(X)$. 
\end{rem}

\vspace{1pc}

\subsection{Construction of a family}

We assume that $\pi:X \to C$ has a section and show that 
$M(\Delta)$ is birational to 
$M(\Delta_0) \times S^nX$,
where $n:=r(\Delta-\Delta_0)$.
Let $\cal E$ be a universal family on 
$M(\Delta_0) \times X$.
Let $(r_1,d_1)$ be the pair of integers such that 
$r_1d-rd_1=-1$ and $0 < r_1<r$,
and let $\cal E_{r_1,d_1}$ be the vector bundle on $X_0 \times_{C_0} X_0$.
Let $j:X_0 \times_{C_0} X_0 \to X_0 \times X$ be the immersion.
We denote the projection $M(\Delta_0) \times X_0 
\to M(\Delta_0)$ by $q_1$
and $M(\Delta_0) \times X_0 \to X_0$ by $q_2$.
By Lemma \ref{lem:elm}, $\cal L:=\Hom_{p_{M(\Delta_0) \times X_0}}
((q_1 \times 1_X)^*\cal E,(q_2 \times 1_X)^*
j_*\cal E_{r_1,d_1})$
is a line bundle on $M(\Delta_0) \times X_0$,
and there is a surjective homomorphism:
$(q_1 \times 1_X)^*\cal E 
 \to (q_2 \times 1_X)^*j_*\cal E_{r_1,d_1}
\otimes p_{M(\Delta_0) \times X_0}^*(\cal L)^{\vee}$. 
Let $p_i:X_0^n:=X_0 \times X_0 \times \dots \times X_0 \to X_0$
be the $i$-th projection, $1 \leq i \leq n$.
Then there is a homomorphism
\begin{equation}
\Lambda:\widetilde{\cal E}  \to 
\oplus_{i=1}^n (q_2 \circ (1_{M(\Delta_0)} \times 
p_i) \times 1_X)^*j_* \cal E_{r_1,d_1}\otimes \cal L_i,
\end{equation}
where $\widetilde{\cal E}$ is the pull-back of $\cal E$
to $M(\Delta_0) \times X_0^n \times X$
and $\cal L_i=(1_{M(\Delta_0)} \times p_i \times 1_X)^*
p_{M(\Delta_0) \times X_0}^*(\cal L)^{\vee}.$ 
We set $\Gamma:=\{(x_1,x_1,\dots,x_n) \in X_0^n|
\pi(x_i)=\pi(x_j) \text{ for some $i \ne j$}\}$.
Then $\Lambda_1:=\Lambda_{|M(\Delta_0) \times 
(X_0^n \setminus \Gamma) \times X}$ is a surjective homomorphism. 
We set $\cal F:=\ker \Lambda_1$.
By Lemma \ref{lem:D}, $\cal F$ is a family of stable vector bundles on $X$.
Hence there is a morphism
$M(\Delta_0) \times (X_0^n \setminus \Gamma)
 \to M(\Delta)$.
By our construction, this morphism is $\frak S_n$-invariant,
and hence we get a morphism 
$\nu:M(\Delta_0) \times (X_0^n /\frak S_n) \to M(\Delta)$.
By our construction, it is injective.
Since $\dim S^n X=2n=\dim M(\Delta)-
\dim M(\Delta_0)$,
ZMT implies that
$M(\Delta_0) \times (X_0^n /\frak S_n)
\to M(\Delta)$ is an immersion.
We set $M(\Delta)^2:=\nu(M(\Delta_0) \times (X_0^n /\frak S_n))$.
We shall show that $M(\Delta)^2$ is dense. 
For this purpose, we shall estimate the dimension of 
$M(\Delta)^1 \setminus M(\Delta)^2$.

\begin{lem}
$\dim(M(\Delta)^1 \setminus M(\Delta)^2) = 2n-1+\dim M(\Delta_0)$.
\end{lem}
 
\begin{pf}
For a $E \in M(\Delta)^1$ and a smooth fibre $l$,
we assume that $E_{|l}$ is not stable.
By the definition of $M(\Delta)^1$, we see that 
$E_{|l} \cong E_1 \oplus E_2$,
where $E_1$ (resp. $E_2$) is a stable vector bundle of rank $r_1$
and degree $d_1$ (resp. rank $r_2$ and degree $d_2$).
We set $E':=\ker(E \to E_1)$.
Then there is an exact sequence 
\begin{equation}\label{eq:1}
0 \to E_1 \to E'_{|l} \to E_2 \to 0.
\end{equation}
Then $E$ is obtained by the inverse transform from $E'$
:
\begin{equation}
0 \to E \to E'(l) \to E_2 \to 0.
\end{equation} 
By \eqref{eq:1},
$E'_{|l}$ is stable or 
$E'_{|l} \cong E_1 \oplus E_2$.
By Lemma \ref{lem:D}, $\Delta(E')=\Delta(E)-1/r$.
Conversely, for $E' \in M(\Delta-1/r)^1$,
we shall consider a surjective homomorphism
$\psi:E' \to F_2$ such that 
the kernel of $E'_{|l} \to F_2$ is stable,
where $F_2$ is a stable vector bundle of rank $r_2$
on a smooth fibre $l$ with degree $d_2$.
If $\ker \psi \otimes \cal O_X(l)$ belongs to 
$M(\Delta)^1 \setminus M(\Delta)^2$,
then (i) $E'_{|l}$ is stable and $E'$ belongs to
$M(\Delta-1/r)^1 \setminus M(\Delta-1/r)^2$,
or (ii) $E'_{|l}$is not stable and $F_2$
is a direct summand of $E'_{|l}$.
Since $\#\{l|\text{$E'_{|l}$ is not stable}\} \leq n-1$,
by using Lemma \ref{lem:elm}, we see that 
\begin{align*}
\dim(M(\Delta)^1 \setminus M(\Delta)^2)& =
\max\{\dim (M(\Delta-1/r)^1 \setminus M(\Delta-1/r)^2)+2,
\dim M(\Delta-1/r)^1+1\}\\
& =2n-1+\dim M(\Delta_0).
\end{align*}
\end{pf}

\begin{thm}\label{thm:1}
$M(\Delta)$ is irreducible
and birational to $M(\Delta_0) \times S^n (J^{d_1}X)$,
where $n:=r(\Delta-\Delta_0)$.
\end{thm}

\begin{pf}
If $\pi:X \to C$ has a section,
we have proved our theorem.
For general cases, we shall consider a Galois covering
$\gamma:C' \to C$ such that $\pi':X \times _C C' \to C'$ has 
a section $\sigma'$.
Let $C_1$ be an open subscheme of $C_0$ such that
$\gamma^{-1}(C_1) \to C_1$ is etale.
We set $X_1':=\pi^{-1}(C_1) \times_{C} C'$.
Let $\cal E_{r_1,d_1}'$ be the vector bundle on 
$X_1' \times_{\gamma^{-1}(C_1)} X_1'$ and 
$j':X_1' \times_{\gamma^{-1}(C_1)} X_1' \cong
X_1' \times_{C_1} X_1 \hookrightarrow X_1' \times X_1$
the inclusion.
Let $X_1' \to J^{d_1}X$ be the morphism induced by $\cal E_{r_1,d_1}'$.
For a $g \in \Gal(C'/C)$, let $\tilde{g}:X_1' \to X_1'$
be the automorphism of $X_1'$ sending
$(x,y) \in \pi^{-1}(C_1) \times_C C'$ to 
$(x+(d_1-1)(\sigma'(g(y))-\sigma'(y)),g(y))$.
Then it defines an action of $\Gal(C'/C)$ to $X_1'$.
By the construction of $\cal E_{r_1,d_1}'$,
we see that $\det(\cal E_{r_1,d_1}')_{|\tilde{g}((x,y))}
\cong \det(\cal E_{r_1,d_1}')_{|(x,y)}$.
Hence $(\cal E_{r_1,d_1}')_{|\tilde{g}((x,y))}
\cong (\cal E_{r_1,d_1}')_{|(x,y)}$.
Thus the morphism $X_1' \to J^{d_1}X$ is $\Gal(C'/C)$-invariant.
Then we get that $X_1'/\Gal(C'/C) \to J^{d_1}X$ is an immersion. 
Replacing $j_* \cal E_{r_1,d_1}$ by    
$j'_* \cal E_{r_1,d_1}'$,
we can construct a family of stable vector bundles $\cal F$
parametrized by $M(\Delta_0) \times 
((X_1')^n \setminus \Gamma')$,
 where $\Gamma'$ is the pull-back of $\Gamma$ to $(X_1')^n$.
Hence we get a morphism 
$M(\Delta_0) \times ((X_1')^n \setminus \Gamma') \to
M(\Delta)$.
By the construction,
$\Gal(C'/C) \times \frak S_n$ acts on 
$((X_1')^n \setminus \Gamma')$, and 
this morphism is $\Gal(C'/C)\times \frak S_n$-invariant.
Hence we get a morphism 
$M(\Delta_0) \times ((J^{d_1}X_1)^n \setminus \Gamma)/\frak S_n \to
M(\Delta)$.
Then we see that $M(\Delta)$
is birationally equivalent to $M(\Delta_0) \times S^n (J^{d_1}X)$.
\end{pf}

\section{Moduli spaces on Del Pezzo surfaces}

\subsection{}
We shall apply Theorem \ref{thm:1} to moduli spaces on Del Pezzo surfaces.
\begin{thm}
We assume that $X=\Bbb P^2$ and set $H:=\cal O_{\Bbb P^2}(1)$.
Then $M_H(r,kH,\Delta)$ is a rational variety if $(r,3k)=1$.
\end{thm}

\begin{pf}
Let $V \subset H^0(\Bbb P^2,\cal O_{\Bbb P^2}(3))$
be a pencil such that every member $D \in V$ is irreducible and 
$\#\{P|P \in \cap_{D \in V} D \}=9$.
Let $\phi:Y \to \Bbb P^2$ be the blow-ups of $\Bbb P^2$
at base points of $V$.
Then there is an elliptic fibration $\pi:Y \to \Bbb P^1$
such that every fibre is isomorphic to a member $D$
of $V$.
We set 
\begin{equation}
N:=\{E \in M_H(r,kH,\Delta)_0| 
\text{$\phi^*E_{|\pi^{-1}(\eta)}$ is stable }
\},
\end{equation}
where $\eta$ is the generic point of $\Bbb P^1$.
Let $E$ be a stable vector bundle of rank $r$ on $\Bbb P^2$
with $c_1(E)=kH$.
Then $\Ext^2(E,E(-3))_0 \cong \Hom(E,E)_0^{\vee}=0$.
Let $D \in V$ be a smooth elliptic curve.
Then we get the surjective homomorphism
$\Ext^1(E,E)_0 \to \Ext^1(E_{|D},E_{|D})_0$.
Hence $\Def(E) \to \Def(E_{|D})$ is submersive.
Since $(r,\deg(E_{|D}))=(r,3k)=1$,
we can deform $E$ to a stable sheaf $F$ such that 
$F_{|D}$ is a stable vector bundle on $D$.
By the openness of stability, $F_{|\pi^{-1}(\eta)}$
is a stable vector bundle.
Hence $N$ is an open dense subscheme of $M_H(r,kH,\Delta)$
and there is an open immersion
$\phi^*:N \to M(r,k\phi^* H,\Delta)$.
By Theorem \ref{thm:1}, $N$ is bitarional to $S^nY$,
where $n=r\Delta-(r^2-1)/2$.
Since $S^nY$ is a rational variety, we get our theorem.
\end{pf} 

\begin{defn}
$Spl(r,c_1,\Delta)$ is the moduli space of simple torsion free sheaves 
$E$ of rank $r$ with $c_1(E)=c_1$ and $\Delta(E)=\Delta$.
\end{defn}
We shall next consider the irreducibility of $Spl(r,c_1,\Delta)$
for Del Pezzo surfaces.
\begin{prop}\label{prop:irr}
Let $\pi:X \to \Bbb P^1$ be a rational elliptic surface
with a section $\sigma$.
For a $c_1 \in \NS(X)$ such that $(c_1,f)$ and r are relatively prime,
we shall consider the moduli space
$M(\Delta)=M(r,c_1,\Delta)$.
Then $M(\Delta)$ is irreducible and rational.
\end{prop}
  
\begin{pf}
We note that $\sigma$ is a $(-1)$-curve.
Let $\phi:X \to Y$ be the contraction of $\sigma$.
Since the characteristic of $\Bbb C$ is 0,
$\pi_* \cal O_X$ is locally free of rank 1, and hence
$\pi_* K_X^{\vee}(\sigma)\cong
\pi_* K_X^{\vee}$.
Then we get that $H^0(Y,K_Y^{\vee}) \cong 
H^0(X,K_X^{\vee}(\sigma))\cong H^0(X,K_X^{\vee})
\cong \Bbb C^{\oplus 2}$.
By the Riemann-Roch theorem, $H^1(Y,K_Y^{\vee})=0$.
Let $\delta:\cal Y \to S$ be a smooth family of 8-points blow-ups
of $\Bbb P^2$ such that $H^1(\cal Y_s,K_{\cal Y_s}^{\vee})=0$
for all $s \in S$ and $\cal Y_{s_0}=Y$ for some $s_0 \in S$.
Let $\xi$ be the generic point of $S$.
By the base change theorem, 
$\delta_*(K_{\cal Y/S}^{\vee})$ is a locally free sheaf of rank
2 and $\delta_*(K_{\cal Y/S}^{\vee})\otimes k(s) \to H^0(K_{Y_s}^{\vee}),
s \in S$ is an isomorphism.
We set $\cal O_{\cal Z}:=\coker(\delta^*
\delta_*(K_{\cal Y/S}^{\vee}) \to K_{\cal Y/S}^{\vee})
\otimes K_{\cal Y/S}$.
Then $\cal O_{\cal Z} \otimes k(s)$ defines a reduced one 
point of $Y_s$.
Thus $\cal Z$ defines a section of $\delta$.
Let $\phi_S:\cal X \to \cal Y$ be the blow-up of $\cal Y$
along $\cal Z$ and set $\epsilon:=\delta \circ \phi$.
Then there is a morphism 
$\pi_S:\cal X \to \Bbb P:=
\Bbb P(\epsilon_*(K_{\cal Y/S}^{\vee}))$,
which defines a family of elliptic fibrations.
Choosing a sufficiently general family,
we may assume that $\pi_{S|\xi}:\cal X_{\xi} \to \Bbb P^1_{k(\xi)}$
is an elliptic surface such that every fibre is irreducible.
Let $\cal O_{\cal X}(1)$ be a relative ample line bundle on $\cal X$ which is 
sufficiently close to the pull-back of an ample line bundle on $\Bbb P$.
For a line bundle $\cal L$ on $\cal X$ such that $c_1(\cal L_{s_0})= c_1$,
we shall consider the relative moduli space
$\cal M(r,\cal L,\Delta) \to S$ of stable sheaves $E$ 
of rank $r$ on $\cal X_s, s \in S$ such that $c_1(E)=\cal L_s$
and $\Delta(E)=\Delta$.
By Maruyama [M1, Cor. 5.9.1, Prop. 6.7], $\cal M(r,\cal L,\Delta)$ is smooth and projective over $S$.
By Theorem \ref{thm:1}, the generic fibre is irreducible,
and hence every fibre is irreducible.
Thus $M(\Delta)$ is irreducible.
Since $M(\Delta)$ contain an irreducible component which
is birational to $S^n X$ for some $n$ ( see the proof of Theorem \ref{thm:1}),
$M(\Delta)$ is a rational variety.
\end{pf}

\begin{lem}\label{lem:spl}
Let $\phi:\widetilde{X} \to X$ be a one point blow-up of
a surface $X$
and $E$ a simple torsion free sheaf of rank $r$ on $X$
which is locally free at the center of the blow-up.
Let $C_1$ be the exeptional divisor of $\phi$
and $\phi^*E \to \cal O_{C_1}^{\oplus k}$, $0<k<r$ a surjective
homomorphism.
We set $E':=\ker(\phi^* E \to \cal O_{C_1}^{\oplus k})$.
Then $E'$ is also a simple torsion free sheaf.
\end{lem}

\begin{pf}
We note that $\Ext^1(\cal O_{C_1}^{\oplus k},E)
\cong H^1(C_1,E^{\vee} \otimes  
\cal O_{C_1}(K_{\widetilde{X}})^{\oplus k})
\cong H^1(C_1,\cal O_{C_1}(-1)^{\oplus rk})=0$.
By the exact sequence
$0 \to E' \to E \to \cal O_{C_1}^{\oplus k} \to 0$,
we see that $\Hom(E,E) \cong \Hom(E',E)$.
Since $\Hom(E',E') \to \Hom(E',E)$ is injective,
we get that $\Hom(E',E')=\Bbb C$.
\end{pf}

\begin{cor}\label{cor:spl}
Let $E$ be a simple torsion free sheaf of rank $r$ on $X$ with 
$c_1(E)=c_1$ and $\Delta(E)=(\Delta)$ which is locally free at the center of a
blow-up $\phi:\widetilde{X} \to X$,
and $E'$ the kernel of a surjective homomorphism
$\phi^* E \to \cal O_{C_1}^{\oplus k}$, $0 \leq k <r$.
We set $\Delta(E')=\Delta'$.
Then,
if $Spl(r,\phi^* c_1-kC_1,\Delta')$ is irreducible,
$Spl(r,c_1,\Delta)$ is also irreducible.
\end{cor}

\begin{pf}
Let $Spl(r,\phi^* c_1,\Delta)^0$ be the open subscheme of 
$Spl(r,\phi^* c_1,\Delta)$ of elements $E$ such that $E_{|C_1} \cong
\cal O_{C_1}^{\oplus r}$.
Then $\phi^*:Spl(r,c_1,\Delta)' \to Spl(r,\phi^* c_1,\Delta)^0$
is an isomorphism, where $Spl(r,c_1,\Delta)'$ is the open dense
subspace of $Spl(r,c_1,\Delta)$ consisting of $E$ such that $E$ is locally free
at the center of the blow-up.
For an $E \in Spl(r,\phi^* c_1,\Delta)^0$,
the quotients $\phi^* E \to \cal O_{C_1}^{\oplus k}$
is parametrized by the Grassmannian variety $G(H^0(C_1,E_{|C_1}),k)$.
Let $Spl(r,\phi^* c_1-kC_1,\Delta')^0$ be the open subscheme of 
$Spl(r,\phi^* c_1-kC_1,\Delta')$ of elements $E'$ 
such that $E'_{C_1} \cong \cal O_{C_1}(1)^{\oplus k} 
\oplus \cal O_{C_1}^{\oplus (r-k)}$.
By using Lemma \ref{lem:spl},
 we can show that there is an open subscheme $U$ of 
$Spl(r,\phi^* c_1-kC_1,\Delta')^0$ and a surjective 
 morphism
$U \to Spl(r,\phi^* c_1,\Delta)^0$
such that every fibre is a Grassmannian variety. 
Hence, the irreducibility of $Spl(r,\phi^* c_1-kC_1,\Delta')$
implies that of $Spl(r,c_1,\Delta)$.
\end{pf}  

\begin{prop}
Let $X$ be a Pel Pezzo surface and
$c_1$ an element of $\NS(X)$.
Then $Spl(r,c_1,\Delta)$ is irreducible.
\end{prop}

\begin{pf}
This follows from Corollary \ref{cor:spl} and Proposition \ref{prop:irr}.
\end{pf}

\section{Moduli spaces on Abelian surfaces}
\subsection{}
For a manifold $V$ and $\alpha \in H^*(V,\Bbb Z)$,
$[\alpha]_i \in H^i(V,\Bbb Z)$ denotes the $i$-th component of $\alpha$.
Let $K(V)$ be the Grothendieck group of $V$.
Let $p:X \to \Spec(\Bbb C)$ be an Abelian surface over $\Bbb C$.
We set
\begin{equation}
\begin{cases}
H^{ev}(X, \Bbb Z):=
H^0(X,\Bbb Z) \oplus  H^2(X,\Bbb Z) \oplus H^4(X,\Bbb Z)\\
H^{odd}(X,\Bbb Z):= H^1(X,\Bbb Z) \oplus  H^3(X,\Bbb Z).
\end{cases}
\end{equation}
Let $E_0$ be an element of $M_H(r,c_1,\Delta)$.
We set 
\begin{equation}
H(r,c_1,\Delta):=
\{\alpha \in H^{ev}(X,\Bbb Z)|[p_*((\chr E_0)\alpha)]_0=0 \}.
\end{equation}
Let $\cal F$ be a quasi-universal family of similitude $\rho$
on $M_H(r,c_1,\Delta) \times X$ [Mu3, Thm. A.5].
Then Mukai [Mu3, Mu5] and Drezet [D, D-N] defines a homomorpism
\begin{equation}
\kappa_2:H(r,c_1,\Delta) 
\to H^2(M_H(r,c_1,\Delta),\Bbb Z)
\end{equation}
such that

\begin{equation}
\kappa_2(\alpha)=\frac{1}{\rho}[p_{M_H(r,c_1,\Delta)*}(\chr(\cal F)\alpha)]_2.
\end{equation}

\begin{rem}
In the notation of Mukai [Mu5, Sect. 5],
$\kappa_2(\alpha)=-\theta_v(\alpha^{\vee})$ and $H(r,c_1,\Delta)=v^{\bot}$,
where $v:=(r,c_1,(c_1^2)/2r-\Delta) \in H^{ev}(X, \Bbb Z)$
is the Chern character of $E_0$.
and $\vee:H^{ev}(X,\Bbb Z) \to 
H^{ev}(X,\Bbb Z)$
is the automorphism sending $\alpha=\alpha_0+\alpha_2+\alpha_4,\;
\alpha_i \in H^{2i}(X,\Bbb Z)$ to
$\alpha^{\vee}=\alpha_0-\alpha_2+\alpha_4$.
Since we used Drezet's notation in [Y2,Y3],
we shall use Drezet's homomorphism in this note.
\end{rem}
We also consider the homomorphism:
\begin{equation}
\kappa_1:H^{odd}(X,\Bbb Z) \to H^1(M_H(r,c_1,\Delta),\Bbb Z)
\end{equation}
such that

\begin{equation}
\kappa_1(\alpha)=\frac{1}{\rho}[p_{M_H(r,c_1,\Delta)*}(\chr(\cal F)\alpha)]_1.
\end{equation}
We note that $\kappa_1$ and $\kappa_2$ do not depend on the choice of $\cal F$.
In this section, we shall prove the following theorem.

\begin{thm}\label{thm:H2}
Let $c_1$ be an element of $\NS(X)$ such that
$c_1 \mod r H^2(X,\Bbb Z)$ is a primitive element of
$H^2(X,\Bbb Z/r \Bbb Z)$ and $H$ a general ample divisor.
We assume that $\dim M_H(r,c_1,\Delta)=2r\Delta+2 \geq 6$.
Let $\frak a:M_H(r,c_1,\Delta) \to \Alb(M_H(r,c_1,\Delta))$ be 
an Albanese map.
Then the following holds.\newline
$(1)$ $\kappa_1$ is an isomorphism and $\kappa_2$ is injective.
\newline
$(2)$
\begin{equation}\label{eq:H2}
\begin{split}
H^2(M_H(r,c_1,\Delta),\Bbb Z) &=
 \kappa_2(H(r,c_1,\Delta))
\oplus \frak a^* H^2(\Alb(M_H(r,c_1,\Delta),\Bbb Z)\\
&=\kappa_2(H(r,c_1,\Delta))\oplus \bigwedge^2 
\kappa_1(H^{odd}(X,\Bbb Z)).
\end{split}
\end{equation}
$(3)$ 
\begin{equation}\label{eq:NS}
\NS(M_H(r,c_1,\Delta)) = 
\kappa_2(H(r,c_1,\Delta)_{alg})
\oplus \frak a^* \NS(\Alb(M_H(r,c_1,\Delta)),
\end{equation}
where $H(r,c_1,\Delta)_{alg}:
=(H^0(X,\Bbb Z) \oplus \NS(X) \oplus H^4(X,\Bbb Z))\cap H(r,c_1,\Delta)$.
\end{thm}
 
\subsection{}
We first assume that $X$ is a product of elliptic curves. 
Let $C_1$ and $C_2$ be elliptic curves and set $X=C_1 \times C_2$.
We set $C_k^i:=C_k$ and $X^i:=C_1^i \times C_2^i$ for $i=0,1,\dots,n,a$,
and $k=1,2$.
Let $\Delta_k^{i,j}$ be the diagonal of $C_k^i \times C_k^j=C_k \times C_k$.
Let $p_k^i$ be a point of $C_k^i$.
We also denote $c_1(\cal O(p_k^i))$ by $p_k^i$.
For simplicity, we denote the pull-backs of $p_k^i$ and $\Delta_k^{i,j}$ to
$X^0 \times Y_0 \times X^a$ by $p_k^i$ and 
$\Delta_k^{i,j}$ respectively.
Let $\Delta_X^{i,j,k}$ be the pull-back of the 
diagonal of $X^i \times X^j \times X^k$ to 
$X^1 \times X^2 \times \dots \times X^n$
and $\Delta_X^{i,j}$ that of $X^i \times X^j$ to 
$X^1 \times X^2 \times \dots \times X^n$.
We set $Z:=\cup_{i<j<k}\Delta_X^{i,j,k}$.
Let $\phi:Y \to (X^1 \times X^2 \times \dots \times X^n)
\setminus Z$
be the blow-up of $(X^1 \times X^2 \times \dots \times X^n)
\setminus Z$
at the subscheme $\cup_{i<j}\Delta_X^{i,j} \setminus Z$,
We set $E^{i,j}:=\phi^{-1}(\Delta_X^{i,j} \setminus Z)$.
For $\alpha \in H^*(X,\Bbb Z)$ and the projection
$\varpi_i:X^0 \times Y_0 \times X^a \to X^i=X$, 
$i=0,1,\dots,n,a$,
we denote the pull-back of $\alpha$ to $X^0 \times Y_0 \times X^a$
by $\alpha^i$.
Then $H^2(\Hilb_X^n,\Bbb Z) \cong H^2(Y,\Bbb Z)^{\frak S_n}$
and $H^2(Y,\Bbb Z)^{\frak S_n}$ is generated by
$\sum_{i=1}^n e^i$, $\sum_{i<j}(f^i \cdot g^j-g^i \cdot f^j)$ 
and $\sum_{i<j}E^{i,j}$
where $e \in H^2(X,\Bbb Z)$ and $f,g \in H^1(X, \Bbb Z)$.
Let $\frak a:X^0 \times \Hilb_X^n \to X^0 \times X$ be the Albanese map
such that $\frak a((x,I_Z))=(x,\sum_{i=1}^n x_i)$
for reduced subscheme $Z=\cup_i\{x_i \}$.

\begin{lem}\label{lem:Ch}
$(1)$ Let $F$ be a vector bundle on $C_2^0 \times C_2^a$ such that
$F_{|\{t\} \times C_2^a}$, $t \in C_2^0$ is a stable vector bundle 
of rank $r$ on $C_2^a$
with $\det F_{|\{t\} \times C_2^a} \cong 
\cal O(\Delta_2^{0,a}+(d-1)p_2^a)_{|\{t\} \times C_2^a}$.
Then,
\begin{equation}
\begin{cases}
c_1(F)=\Delta_2^{0,a}+(d-1)p_2^a+(r_1-1+kr)p_2^0,\; k \in \Bbb Z\\
\chr_2(F)=\frac{1}{2r}(c_1(F)^2).
\end{cases}
\end{equation}
If $k=0$, then $\chr_2(F)=d_1 p_2^0 \cdot p_2^a$.\newline
$(2)$ Let $F_i$ ($1 \leq i \leq n$) be a vector bundle on 
$C_2^i \times C_2^a$ such that
$F_{i|\{t\} \times C_2^a}$, $t \in C_2^i$ is a stable vector bundle of rank $r_2$ on $C_2^a$
with $\det F_{i|\{t\} \times C_2^a} \cong 
\cal O(\Delta_2^{i,a}+(d_2-1)p_2^a)_{|\{t\} \times C_2^a}$.
Then, 
\begin{equation}
\begin{cases}
c_1(F_i)=\Delta_2^{i,a}+(d_2-1)p_2^a+(r_1-1+kr)p_2^i,\; k \in \Bbb Z\\
\chr_2(F_i)=\frac{1}{2r_2}(c_1(F_i)^2).
\end{cases}
\end{equation}
If $k=0$, then $\chr_2(F_i)=d_1 p_2^i \cdot p_2^a$.
\end{lem}

\begin{pf}
We shall only prove (1).
We set $c_1(F)=\Delta_2^{0,a}+(d-1)p_2^a+(r_1-1+x)p_2^0$, $x \in \Bbb Z$.
Since $F_{|\{t\} \times C_2^a}$, $t \in C_2^0$ is a stable vector bundle,
$\Delta(F)=c_2(F)-(c_1(F)^2)(r-1)/2r=0$.
Hence we get that $\chr_2(F)=-(c_2(F)-(c_1(F)^2)/2)=(c_1(F)^2)/2r$.
We note that $c_2(F)=(d(r_1+x)-1)(r-1)/r$ is an integer.
Hence $d(r_1+x)-1=rd_1+rx$ is a multiple of $r$.
Since $(r,d)=1$, $x$ is a multiple of $r$.
We also see that $(c_1(F)^2)/2r=d_1 p_2^0 \cdot p_2^a$ for the case $x=0$.
\end{pf}

Let $F$ and $F_i$ be vector bundles in Lemma \ref{lem:Ch}
and assume that $k=0$.
We also denote the pull-backs of $F$ and $F_i$ to 
$C_2^0 \times C_2^i \times C_2^a$
by $F$ and $F_i$ respectively.
Let $q_{C_2^0 \times C_2^i}:
C_2^0 \times C_2^i \times C_2^a \to C_2^0 \times C_2^i$
be the projection.
We set $\cal L:=\Hom_{q_{C_2^0 \times C_2^a}}(F,F_i)$.
Then $c_1(\cal L)=-\Delta_2^{0,i}$.
\begin{pf}
By using the Grothendieck-Riemann-Roch theorem and the above lemma,
we see that
\begin{align*}
c_1(\cal L)&=[q_{C_2^0 \times C_2^i*}(\chr(F^{\vee})\chr( F_i))]_2\\
&=[q_{C_2^0 \times C_2^i*}(r-c_1(F)+\frac{1}{2r}(c_1(F)^2))
(r_2-c_1(F_i)+\frac{1}{2r_2}(c_1(F_i)^2))]_2\\
&=[q_{C_2^0 \times C_2^i*}(rr_2+(rc_1(F_i)-r_2c_1(F))+
\frac{1}{2rr_2}((rc_1(F_i)-r_2c_1(F))^2))]_2\\
&=\frac{1}{2rr_2}[q_{C_2^0 \times C_2^i*}((rc_1(F_i)-r_2c_1(F))^2)]_2\\
&=-\Delta_2^{0,i}.
\end{align*}
\end{pf}
Let $Y_0$ be the complement of the closed subset
$W:=\cup_{i<j<k}(\widetilde{\Delta}_1^{i,j} \cap 
\widetilde{\Delta}_1^{j,k}) \cup \cup_{i<j}
(\widetilde{\Delta}_1^{i,j} \cap E^{i,j})$ of $Y$,
where $\Delta_1^{i,j}=\widetilde{\Delta}_1^{i,j} \cup E^{i,j}$.
Since $\codim W=2$, 
$H^2(X^0 \times Y_0,\Bbb Z) \cong H^2(X^0 \times Y,\Bbb Z)$.

We shall construct a family of stable sheaves on $X$ parametrized by 
$X^0 \times Y_0$.
For simplicity, we denote the pull-backs of $F$ and $F_i$
to $X^0 \times Y_0 \times X^a$ by $F$ and $F_i$ respectively.
Then there is a homomorphism:

\begin{equation}
\Lambda:F \otimes \cal O(\Delta_1^{0,a}-p_1^a) \to 
\oplus_{i=1}^n (F_{i|\Delta_1^{i,a}} \otimes L^i),
\end{equation}
where $L^i$ is a line bundle on $X^0 \times Y_0 \times X^a$
such that $c_1(L^i)=\Delta_1^{0,i}-p_1^i+\Delta_2^{0,i}$.
Let $\cal E$ be the kernel of this homomorphism and $\cal Q$
the cokernel. 
Then $\cal Q \cong 
\oplus_{i<j}((F_i/G_j)_{|\Delta_1^{i,a} \cap \widetilde{\Delta}_1^{i,j}}
\otimes L^i
\oplus(F_i \otimes L^i_{|\Delta_1^{i,a}} 
\otimes \cal O_{E^{i,j}}))$,
where $G_i:=\ker(F \to F_i)$.
We first assume that $r_1 \leq r_2$.
Then $G_{j|\Delta_1^{i,a}} \to F_{i|\Delta_1^{i,a}}$
 is injective and $(F_i/G_j)_{|\Delta_1^{i,a}}$
is flat over $X^0 \times Y_0$. 
Hence we see that 
\begin{equation}
\Tor^{\cal O_{X^0 \times Y_0}}_2((F_i/G_j)_{|\Delta_1^{i,a} 
\cap \widetilde{\Delta}_1^{i,j}},k(x))=0,\;
x \in X^0 \times Y_0.
\end{equation}
Since $F_i \otimes \cal O(\Delta_2^{0,i})_{|\Delta_1^{i,a}}$ is also flat 
over $X^0 \times Y_0$,
we get that
\begin{equation}
 \Tor^{\cal O_{X^0 \times Y_0}}_2
(F_i \otimes \cal O(\Delta_2^{0,i})_{|\Delta_1^{i,a}} 
\otimes \cal O_{E^{i,j}},k(x))=0,\;
x \in X^0 \times Y_0.
\end{equation}
Hence we see that
$\Tor^{\cal O_{X^0 \times Y_0}}_1(\im(\Lambda),k(x))=0$, which implies that
$\cal E$ is flat over $X^0 \times Y_0$ and 
$\cal E \otimes k(x)$ is torsion free.
Then $\cal E$ defines a family of stable sheaves on $X$ parametrized by 
$ X^0 \times Y_0.$
It defines a morphism $X^0 \times Y_0 \to M(r,c_1,\Delta)$,
which is $\frak S_n$-invariant.
Hence we get a morphism $\nu:X^0 \times (Y_0/\frak S_n) \to
M(r,c_1,\Delta)$. 

Let $\overline{\kappa_2}:H(r,c_1,\Delta) \to 
H^2(X^0 \times Y_0,\Bbb Z)/\frak a^* 
H^2(\Alb(X^0 \times \Hilb_X^n),\Bbb Z)$
be the homomorphism sending $\alpha \in H(r,c_1,\Delta)$
to $[p_{X^0 \times Y_0*} (\chr(\cal E) \alpha)]_2 \mod \frak a^* 
H^2(\Alb(X^0 \times \Hilb_X^n),\Bbb Z)$.
Since $\kappa_2$ does not depend on the choice of quasi-universal
families, we shall compute the image of $\overline{\kappa_2}$.

\begin{align*}
\chr(\cal E) 
&= \chr(F \otimes \cal O(\Delta_1^{0,a}-p_1^a))-
\sum_{i=1}^n \chr(F_i \otimes L^i_{|\Delta_1^{i,a}})+
\sum_{i<j}\chr(F_i\otimes L^i_{|\Delta_1^{i,a}} 
\otimes \cal O_{E^{i,j}}) \\
&\phantom{ (r+c_1(F)+d_1p_2^0\cdot p_2^a)(1+\Delta_1^{0,a})-}
+\sum_{i<j}\chr(F_i/G_j \otimes L^i_{|\Delta_1^{i,a}}
 \otimes \cal O_{\widetilde{\Delta}_1^{i,j}})\\
&= (r+c_1(F)+d_1p_2^0\cdot p_2^a)(1+\Delta_1^{0,a}-p_1^a-p_1^0 \cdot p_1^a)
-\sum_{i=1}^n \Delta_1^{i,a}
(r_2+c_1(F_i)+d_1p_2^i\cdot p_2^a)(\chr L^i)\\
&\phantom{ (r+c_1(F)}+
\sum_{i<j}\chr(F_i \otimes L^i_{|\Delta_1^{i,a}} 
\otimes \cal O_{E^{i,j}})+\sum_{i<j}\chr(F_i/G_j 
\otimes L^i_{|\Delta_1^{i,a}}
 \otimes \cal O_{\widetilde{\Delta}_1^{i,j}}).
\end{align*}
Since $[p_{X^0 \times Y_0*} 
(\chr(F \otimes \cal O(\Delta_1^{0,a}-p_1^a)) \alpha^a)]_2
\equiv 0,\;
\sum_{i=1}^n\Delta_1^{0,i}-p_1^i \equiv 0
 \mod \frak a^* H^2(\Alb(X^0 \times \Hilb_X^n),\Bbb Z)$, we get that
\begin{align*}
\overline{\kappa_2}(\alpha)= & -\sum_{i=1}^n
 [p_{X^0 \times Y_0*} (\Delta_1^{i,a}
(r_2+c_1(F_i)+d_1p_2^i\cdot p_2^a)(1+\Delta_2^{0,i})\alpha^a)]_2\\
&\phantom{-\sum_{i=1}^n [p_{X^0 \times Y_0*}}
 +\sum_{i<j}[p_{X^0 \times Y_0*}(\chr(F_i \otimes 
\cal O(\Delta_2^{0,i})_{|\Delta_1^{i,a}} 
\otimes \cal O_{E^{i,j}})\alpha^a)]_2\\
&\phantom{-\sum_{i=1}^n [p_{X^0 \times Y_0*}(\Delta_1^{i,a}}
+\sum_{i<j}[p_{X^0 \times Y_0*}(\chr(F_i/G_j 
\otimes \cal O(\Delta_2^{0,i})_{|\Delta_1^{i,a}}
 \otimes \cal O_{\widetilde{\Delta}_1^{i,j}})\alpha^a)]_2. 
\end{align*}
Let $\alpha=x_1+x_2p_1+x_3p_2+x_4p_1 \cdot p_2
+D$ be an element of $H(r,c_1,\Delta)$,
 $D \in H^1(C_1, \Bbb Z) \otimes H^1(C_2, \Bbb Z)$.
Then we see that
$0=[p_*((\chr E_0) \alpha)]_0=
[p_*((r+dp_2-r_2np_1-d_2n p_1 \cdot p_2)\alpha)]_0
=-d_2n x_1-r_2 n x_3+d x_2+rx_4$.
Thus $\alpha$ satisfies 
\begin{equation}\label{eq:perp}
d x_2+r x_4=d_2n x_1+r_2 n x_3.
\end{equation}
By a simple calculation, we get that
\begin{equation}
\left\{
\begin{aligned}
&[p_{X^0 \times Y_0*}(\chr (F_i \otimes 
\cal O(\Delta_2^{0,i})_{|\Delta_1^{i,a}})]_2 = d_2 \Delta_2^{0,i}+d_1 p_2^i\\
&[p_{X^0 \times Y_0*}(\chr (F_i \otimes 
\cal O(\Delta_2^{0,i})_{|\Delta_1^{i,a}})p_2^a )]_2 
=r_2 \Delta_2^{0,i}+r_1 p_2^i\\
&[p_{X^0 \times Y_0*}(\chr (F_i \otimes 
\cal O(\Delta_2^{0,i})_{|\Delta_1^{i,a}})p_1^a )]_2  =d_2 p_1^i\\
&[p_{X^0 \times Y_0*}(\chr (F_i \otimes 
\cal O(\Delta_2^{0,i})_{|\Delta_1^{i,a}})D^a )]_2  =D^i\\
&[p_{X^0 \times Y_0*}(\chr (F_i \otimes 
\cal O(\Delta_2^{0,i})_{|\Delta_1^{i,a}})(p_1^a \cdot p_2^a))]_2 =r_2 p_1^i,
\end{aligned}
\right.
\end{equation}

\begin{equation}
\left\{
\begin{aligned}
&[p_{X^0 \times Y_0*}(\chr(F_i/G_j \otimes 
\cal O(\Delta_2^{0,i})_{|\Delta_1^{i,a}} 
\otimes \cal O_{\widetilde{\Delta}_1^{i,j}}))]_2=
(2d_2-d)\widetilde{\Delta}_1^{i,j} \\
&[p_{X^0 \times Y_0*}(\chr(F_i/G_j \otimes 
\cal O(\Delta_2^{0,i})_{|\Delta_1^{i,a}} 
 \otimes \cal O_{\widetilde{\Delta}_1^{i,j}})p_2^a)]_2
=(2r_2-r) \widetilde{\Delta}_1^{i,j}\\
&[p_{X^0 \times Y_0*}(\chr(F_i/G_j \otimes 
\cal O(\Delta_2^{0,i})_{|\Delta_1^{i,a}}  
\otimes \cal O_{\widetilde{\Delta}_1^{i,j}})p_1^a)]_2=0\\
&[p_{X^0 \times Y_0*}(\chr(F_i/G_j
 \otimes \cal O(\Delta_2^{0,i})_{|\Delta_1^{i,a}}  
\otimes \cal O_{\widetilde{\Delta}_1^{i,j}})D^a)]_2=0\\
&[p_{X^0 \times Y_0*}(\chr(F_i/G_j \otimes 
\cal O(\Delta_2^{0,i})_{|\Delta_1^{i,a}}  
\otimes \cal O_{\widetilde{\Delta}_1^{i,j}})(p_1^a \cdot p_2^a))]_2=0,
\end{aligned}
\right.
\end{equation}
and
\begin{equation}
\left\{
\begin{aligned}
&[p_{X^0 \times Y_0*}(\chr(F_i \otimes \cal O(\Delta_2^{0,i})_{|\Delta_1^{i,a}} 
\otimes \cal O_{E^{i,j}}))]_2=d_2 E^{i,j}\\
&[p_{X^0 \times Y_0*}(\chr(F_i \otimes \cal O(\Delta_2^{0,i})_{|\Delta_1^{i,a}} 
 \otimes \cal O_{E^{i,j}})p_2^a)]_2=r_2 E^{i,j}\\
&[p_{X^0 \times Y_0*}(\chr(F_i \otimes \cal O(\Delta_2^{0,i})_{|\Delta_1^{i,a}}  
\otimes \cal O_{E^{i,j}})p_1^a)]_2=0\\
&[p_{X^0 \times Y_0*}(\chr(F_i \otimes \cal O(\Delta_2^{0,i})_{|\Delta_1^{i,a}}  
\otimes \cal O_{E^{i,j}})D^a)]_2=0\\
&[p_{X^0 \times Y_0*}(\chr(F_i \otimes \cal O(\Delta_2^{0,i})_{|\Delta_1^{i,a}}  
\otimes \cal O_{E^{i,j}})(p_1^a \cdot p_2^a))]_2=0,
\end{aligned}
\right.
\end{equation}
where $D \in H^1(C_1, \Bbb Z) \otimes H^1(C_2, \Bbb Z)$.
Hence we get that
\begin{equation}
\begin{align*}
\overline{\kappa_2}(\alpha)=&
-\sum_{i=1}^n(d_2 x_1+r_2 x_3)\Delta_2^{0,i}-
\sum_{i=1}^n(d_2 x_2+r_2 x_4)p_1^i
-\sum_{i=1}^n(d_1 x_1+r_1 x_3)p_2^i-\sum_{i=1}^n D^i\\
&+\sum_{i<j}((2d_2-d)x_1+(2r_2-r) x_3)\widetilde{\Delta}_1^{i,j}
+\sum_{i<j}(d_2 x_1+r_2 x_3)E.
\end{align*}
\end{equation}
We note that
\begin{equation}
\left\{
\begin{aligned}
&\sum_{i=1}^n \Delta_2^{0,i} \equiv \sum_{i=1}^n p_2^i 
\mod \frak a^* H^2(\Alb(X^0 \times \Hilb_X^n),\Bbb Z)\\
&\sum_{i<j}\Delta_1^{i,j} \equiv 2(n-1)\sum_{i=1}^n p_1^i 
\mod \frak a^* H^2(\Alb(X^0 \times \Hilb_X^n),\Bbb Z)\\
& \widetilde{\Delta}_1^{i,j}=\Delta_1^{i,j}-E^{i,j}.
\end{aligned}
\right.
\end{equation}
Therefore we get that
\begin{equation}
\overline{\kappa_2}(\alpha)
=y_1(\sum_{i=1}^n p_2^i) +y_2 (\sum_{i=1}^n p_1^i)
+y_3(\sum_{i<j} E^{i,j})-\sum_{i=1}^n D^i,
\end{equation}
where
\begin{equation}
\begin{cases}
y_1=-(dx_1+r x_3)\\
y_2=-\{(d_2x_2+r_2 x_4)-2(n-1)((2d_2-d)x_1+(2r_2-r)x_3)\}\\
y_3=(d_1x_1+r_1x_3)\\
y_4=dx_2+rx_4-n(d_2x_1+r_2x_3).
\end{cases}
\end{equation}
Since $dr_1-rd_1=d_2r-dr_2=1$,
the homomorphism $\psi:\Bbb Z^{\oplus 4} \to \Bbb Z^{\oplus 4}$
sending $(x_1,x_2,x_3,x_4)$ to
$(y_1,y_2,y_3,y_4)$ is an isomorphism.
The condition \eqref{eq:perp} implies that $y_4=0$.
Therefore,
\begin{equation}
\overline{\kappa_2}:
H(r,c_1,\Delta) \to 
H^2(X^0 \times Y_0,\Bbb Z)^{\frak S_n}/
\frak a^* H^2(\Alb(X^0 \times \Hilb_X^n),\Bbb Z)
\end{equation}
is an isomorphism.
Since $H^2(X^0 \times Y_0,\Bbb Z)^{\frak S_n} \cong
H^2(X^0 \times \Hilb_X^n,\Bbb Z)$,
we get that 
\begin{equation}
H(r,c_1,\Delta) \to 
H^2(X^0 \times \Hilb_X^n,\Bbb Z)/
\frak a^* H^2(\Alb(X^0 \times \Hilb_X^n),\Bbb Z)
\end{equation}
is an isomorphism.

We next treat the case $r_1>r_2$.
Since $G_j \to F_i$ is surjective,
we get that
\begin{equation}
\begin{align*}
\overline{\kappa_2}(\alpha)=&
-\sum_{i=1}^n(d_2 x_1+r_2 x_3)\Delta_2^{0,i}-
\sum_{i=1}^n(d_2 x_2+r_2 x_4)p_1^i
-\sum_{i=1}^n(d_1 x_1+r_1 x_3)p_2^i\\
&-\sum_{i=1}^n D^i+\sum_{i<j}(d_2 x_1+r_2 x_3)E.
\end{align*}
\end{equation}
In the same way as in the case $r_1 \leq r_2$,
we see that
\begin{equation}
H(r,c_1,\Delta) \to 
H^2(X^0 \times \Hilb_X^n,\Bbb Z)/
\frak a^* H^2(\Alb(X^0 \times \Hilb_X^n),\Bbb Z)
\end{equation}
is an isomorphism.

Therefore $\kappa_2$ is injective and $H^2(M_H(r,c_1,\Delta))$ is
generated by $\im(\kappa_2)$ and $\im(\frak a)$.
By using similar computations,
we see that $\kappa_1$ is an isomorphism.
Hence Theorem \ref{thm:H2} (1), (2) hold for this case.

\subsection{}
We next treat general cases.
Replacing $c_1$ by $c_1+r c_1(H)$, we may assume that
$c_1$ belongs to the ample cone.

\begin{prop}\label{prop:1}
Let $(X,L)$ be a pair consisting of abelian surface $X$ and an ample divisor
$L$ of type $(d_1,d_2)$, where $d_1$ and $d_2$ are positive integers of
$d_1|d_2$ and $(r, d_1)=1$.
Then Theorem \ref{thm:H2} $(1)$, $(2)$ hold for $M_H(r,c_1(L),\Delta)$,
where $H$ is a general polarization. 
\end{prop}
\begin{pf}
Let $(X,L)$ be a pair consisting of abelian surface $X$ and an ample divisor
$L$ of type $(d_1,d_2)$, where $d_1$ and $d_2$ are positive integers of
$d_1|d_2$ and $(r, d_1)=1$.
We shall choose an ample line bundle $H$ on $X$ 
which is not lie on walls.
Let $T$ be a connected smooth curve and
$(\cal X, \cal L)$ a pair of a smooth family of
abelian surface $p_T:\cal X \to T$ and a relatively ample line bundle
$\cal L$ of type $(d_1,d_2)$.
For points $t_0, t_1 \in T$, we assume that $(\cal X_{t_0},\cal L_{t_0})
=(X,L)$ and $\cal X_{t_1}$ is an abelian surface of $\NS(\cal X_{t_1})
\cong \Bbb Z$.
Let $g:Pic_{\cal X/T} \to T$ be the relative Picard scheme.
We denote the connected component of $Pic_{\cal X/T}$ containing
the section of $g$ which corresponds to the family $\cal L$ by 
$Pic^{\xi}_{\cal X/T}$.
Since $Pic^0_{\cal X/T} \cong Pic^{\xi}_{\cal X/T}$,
$Pic^{\xi}_{\cal X/T} \to T$ is a smooth morphism.
Let $h:\overline{\cal M}_{\cal X/T}(r,\xi,\Delta) \to T$
be the moduli scheme parametrizing $S$-equivalence classes
of $\cal L_t$-semi-stable sheaves $E$ on $\cal X_t$ with
$(\rk(E),c_1(E),\Delta(E))=(r,c_1(\cal L_t),\Delta)$ [Ma1].
Let $D$ be the closed subset of $\overline{\cal M}_{\cal X/T}(r,\xi,\Delta)$
consisting of properly $\cal L_t$-semi-stable sheaves on $\cal X_t$.
Since $h$ is a proper morphism, $h(D)$ is a closed subset of $T$.
Since $h(D)$ does not contain $t_1$ and $T$ is an irreducible curve,
$h(D)$ is a finite point set.
Replacing $T$ by the open subscheme $T \setminus (h(D)\setminus \{t_0\})$,
we may assume that $\cal L_t$-semi-stable sheaves are $\cal L_t$-stable
for $t \ne t_0$.  
Let $s:\cal Spl_{\cal X/T}(r,\xi,\Delta) \to T$ be the moduli of 
simple sheaves $E$ on $\cal X_t, t \in T$ with
$(\rk(E),c_1(E),\Delta(E))=(r,c_1(\cal L_t),\Delta)$ [A-K, Thm. 7.4].
Let $U_1$ be the closed subset of $s^{-1}(T \setminus \{t_0 \})$
consisting of simple sheaves on $\cal X_t$, $t \in T \setminus \{t_0 \}$ 
which are not stable with respect to $\cal L_t$ 
and $\overline{U_1}$
the closure of $U_1$ in $\cal Spl_{\cal X/T}(r,\xi,\Delta)$.
Let $U_2$ be the closed subset of $s^{-1}(t_0)$
consisting of simple sheaves which are not semi-stable with respect to $H$.
Then we can show that $\overline{U_1} \cap s^{-1}(t_0)$ is a subset of
$U_2$ (see the second paragraph of the proof of Lemma \ref{lem:1}).
We set $\cal M:=\cal Spl_{\cal X/T}(r,\xi,\Delta) 
\setminus (\overline{U_1} \cup U_2)$.
Then $\cal M$ is an open subspace of $\cal Spl_{\cal X/T}(r,\xi,\Delta)$
which is of finite type and contains all $H$-stable sheaves on $\cal X_{t_0}$. 
By using valuative criterion of separatedness and properness,
we get that $s:\cal M \to T$ is a proper morphism.
In fact, since $\cal M \times_T (T \setminus \{t_0\}) \to T \setminus \{t_0\}$
is proper, it is sufficient to check these properties near the fibre 
$\cal X_{t_0}$.   
The separatedness follows from base change theorem and stability
with respect to $H$ (cf. [A-K, Lem. 7.8]),
and the properness follows from the following lemma (Lemma \ref{lem:1})
and the projectivity of $\cal M_{t_0}$.
Since $Pic^{\xi}_{\cal X/T} \to T$ is a smooth morphism,
[Mu2, Thm. 1.17] implies that $s:\cal M \to T$
is a smooth morphism.
Let $\frak a_T:\cal M \to \Alb_{\cal M/T}$ be the family of Albanese map
over $T$.
Let $\cal F_T$ be a quasi-universal family of similitude $\rho$
on $\cal M \times_T \cal X$
and we shall consider the homomorphism
\begin{equation}
\begin{cases}
\kappa_{1,t}:H^{odd}(\cal X_t,\Bbb Z) \to H^1(\cal M_t,\Bbb Z)\\
\kappa_{2,t}:H(r,c_1(\cal L_t),\Delta) \to H^2(\cal M_t,\Bbb Z)
\end{cases}
\end{equation}
such that
$\kappa_{i,t}(\alpha_{i,t})=
\frac{1}{\rho}[p_{\cal M_t*}((\chr \cal F_t)\alpha_t)]_i
$, where $\alpha_{1,t} \in H^{odd}(\cal X_t,\Bbb Z)$,
$\alpha_{2,t} \in H(r,c_1(\cal L_t),\Delta)$.
We assume that $\cal X_{t_0}$ is a product of elliptic curves.
Since $p_T$ and $s$ are smooth, Theorem \ref{thm:H2} (1),(2) for the pair 
$(\cal X_{t_0},\cal L_{t_0})$ imply
 that Theorem \ref{thm:H2} (1),(2) also hold for
all pairs $(\cal X_t, \cal L_t)$, $t \in T$.
By the connectedness of the moduli of $(d_1,d_2)$-polarized abelian surfaces
(cf. [L-B, 8]),
\eqref{eq:H2} holds for all pairs $(X,L)$ of 
$(d_1,d_2)$-polarized abelian surfaces.
\end{pf}
The following is due to Langton [L].
 
\begin{lem}\label{lem:1}
Let $R$ be a discrete valuation ring,
$K$ the quotient field of $R$, and $k$ the residue field of $R$.
Let $\Spec(R) \to T$ be a dominant morphism such that 
$\Spec(k) \to T$ defines the point $t_0$.
For a stable sheaf $E_K$ on $X_K$,
there is a $R$-flat coherent sheaf $E$ on $X_R$ such that $E \otimes_R K=E_K$
and $E \otimes_R k$ is a $H$-stable sheaf.
\end{lem}
\begin{pf}
Let $E^0$ be an $R$-flat coherent sheaf on $X_R$
such that $E^0 \otimes_R K=E_K$ and $E^0_k:=E^0 \otimes_R k$ is torsion free.
If $E_k^0$ is $H$-stable, then we put $E=E^0$.
We assuime that $E^0_k$ is not $H$-stable.
Let $F ^0_k (\subset E^0_k)$ be the first filter of the Harder-Narasimhan filtration 
of $E^0_k$ with respect to $H$.
We set $E^1:=\ker(E^0 \to E^0_k/F^0_k)$.
Then $E^1$ is an $R$-flat coherent sheaf on $X_R$ with $E^1_K=E_K$.
If $E^1_k$ is not $H$-stable,
then we shall consider the first filter $F_k^1$ of the
Harder-Narasimhan filtration of $E^1_k$   
and set $E^2:=\ker(E^1 \to E^1_k/F^1_k)$.
Continuing this procedure successively,
we obtain a decreasing sequence of $R$-flat coherent sheaves on $X_R$:
$E^0 \supset E^1 \supset E^2 \supset \cdots$.
We assume that this sequence is infinite.
Then in the same way as in [L, Lem. 2],
we see that there is an integer $i$ such that $E^i \otimes_R \widehat{R}$
has a subsheaf $F$ of rank $r'$
with $F \otimes_R k=F_k^i$, where $\widehat{R}$ is the completion of $R$.

We set $\widehat{K}:=K \otimes_R \widehat{R}$
and $D:=\det(E^i\otimes_R \widehat{R})^{\otimes r'}\otimes
 \det(F)^{\otimes (-r)}$.
Let $P(x)$ be the Hilbert polynomial of $D$ with respect to 
$\cal L_{\widehat{R}}$.
Let $V$ be a locally free sheaf on $\cal X$ such that 
there is a surjective homomorphism $V \otimes_{\cal O_T} \widehat{R} \to D$,
and we shall consider the quot scheme
$\cal Q:=\Quot_{V/\cal X/T}^{P(x)}$.
Then $D$ defines a morphism $\tau:\Spec(\widehat{R}) \to \cal Q$
such that $D=(\tau \times_T 1_{\cal X})^* \cal D$,
where $\cal D$ is the universal quotient.  
Let $\cal Q_0$ be the connected component of $\cal Q$
which contains the image of $\Spec(\widehat{R})$.
Since $\Spec(\widehat{R}) \to T$ is dominant,
$\frak q:\cal Q_0 \to T$ is dominant, and hence surjective.
Since $E^i_{\widehat{K}} \cong E_K \otimes_K \widehat{K}$
is a stable sheaf on $X_{\widehat{K}}$,
$(\cal D_{q_1},\cal L_{q_1})=(\cal D_{\widehat{K}},\cal L_{\widehat{K}})>0$,
where $q_1$ is a point of $\frak q^{-1}(t_1)$.
Since $\NS(\cal X_{t_1}) \cong \Bbb Z$, we get that
$c_1(\cal D_{q_1})=l c_1(\cal L_{q_1})$, $l>0$.
Hence we obtain that $(\cal D_{\tau(t_0)}^2)>0$ and 
$(\cal D_{\tau(t_0)},\cal L_{\tau(t_0)})>0$.
By the Riemann-Roch theorem and the Serre duality,
we see that $\cal D_{\tau(t_0)}$ is an effective divisor.
Therefore $(\cal D_{\tau(t_0)},H)>0$,
which is a contradiction.
Hence there is an integer $n$ such that $E^n \otimes_R k$
is $H$-stable.
\end{pf}
{\it Proof of Theorem} \ref{thm:H2} (3).
Let $\kappa_2':H(r,c_1,\Delta)\otimes \Bbb C \to H^2(M(r,c_1,\Delta),\Bbb C)$
be the homomorphism induced by $\kappa_2$.
We note that $H^{2,0}(X)$ and $H^{0,2}(X)$ are subsets of 
$H(r,c_1,\Delta)\otimes \Bbb C$.   
Since $\chr_i(\cal F)$ is of type $(i,i)$, we see that
\begin{equation}
\begin{cases}
\kappa_2'(H^{2,0}(X)) \subset H^{2,0}(M_H(r,c_1,\Delta))\\
\kappa_2'(\oplus_{p=0}^2 H^{p,p}(X)) \subset H^{1,1}(M_H(r,c_1,\Delta))\\
\kappa_2'(H^{0,2}(X)) \subset H^{0,2}(M_H(r,c_1,\Delta)).
\end{cases}
\end{equation}
Since $H(r,c_1,\Delta)\otimes \Bbb C=H^{2,0}(X)\oplus 
(\oplus_{p=0}^2 H^{p,p}(X)) \cap H(r,c_1,\Delta)\otimes \Bbb C \oplus
H^{0,2}(X)$ and $\frak a^*$ preserves the type,
we obtain that
$$
H^{1,1}(M_H(r,c_1,\Delta))=\kappa_2'((\oplus_{p=0}^2 H^{p,p}(X))\cap 
H(r,c_1,\Delta)\otimes \Bbb C)\oplus 
\frak a^* (H^{1,1}(\Alb(M_H(r,c_1,\Delta)).
$$
Hence we get Theorem \ref{thm:H2} (3).
\qed

Combining [Y4, Thm. 2.1] with the proof of Proposition \ref{prop:1},
we get the following theorem.
\begin{thm}\label{thm:B}
Let $X$ be an abelian surface defined over $\Bbb C$ and
$c_1 \in \NS(X)$ a primitive element.
Then 
$$
P(M_H(2,c_1,\Delta),z)=P(M_H(1,0,2\Delta),z)
$$ 
for a general polarization $H$,
where $P(\quad \; ,z)$ is the Poincar\'{e} polynomial.
\end{thm}

\subsection{}
We shall next consider the Albanese variety of $M_H(r,c_1,\Delta)$.
Let $\cal P$ be the Poincar\'{e} line bundle on
$\widehat{X} \times X$,
where $\widehat{X}$ is the dual of $X$.
For an element $E_0 \in M_H(r,c_1,\Delta)$,
let $\alpha_{E_0}:M_H(r,c_1,\Delta) \to X$ be the morphism
sending $E \in M_H(r,c_1,\Delta)$ to 
$\det p_{\widehat{X}!}((E-E_0)\otimes(\cal P-\cal O_{\widehat{X} \times X}))
\in \Pic^0(\widehat{X})=X$,
and $\det_{E_0}:M_H(r,c_1,\Delta) \to \widehat{X}$
the morphism sending $E$ to $\det E \otimes \det E_0^{\vee}
\in \widehat{X}$ (cf. [Y3, Sect. 5]).
We shall show that $\frak a_{E_0}:=\det_{E_0} \times \alpha_{E_0} $
is the Albanese map of $M_H(r,c_1,\Delta)$.
Let $B$ be an effective divisor on $X$.
Then we see that 
\begin{align*}
&\det p_{\widehat{X}!}((E-E_0)\otimes \cal O_B \otimes 
(\cal P-\cal O_{\widehat{X} \times X}))\\
=&
\det p_{\widehat{X}!}((\det E_{|B}-\det E_{0|B})\otimes 
(\cal P-\cal O_{\widehat{X} \times X}))\\
=& \zeta (\det\nolimits_{E_0}(E)),
\end{align*}
where $\zeta:\widehat{X} \to X$ is the morphism sending $L \in \widehat{X}$
to $\otimes_i \cal P_{\widehat{X} \times \{x_i\}} \in \Pic^0(\widehat{X})=X$, 
$L \cdot B = \sum_i x_i$.
Therefore if $\frak a_{E_0}$ is the Albanese map for $M_H(r,c_1,\Delta)$,
then $\frak a_{E_0(B)}$ is the Albanese map for $M_H(r,c_1+rc_1(\cal O_X(B)),
\Delta)$.
Hence we may assume that $c_1$ belongs to the ample cone.
In the notation of Proposition \ref{prop:1},
we assume that there is a section $\sigma:T \to \cal M$ of $s$.
Then we can also construct a morphism 
$\frak a_{\sigma}:\cal M \to \Pic^0_{\cal X/T} \times_T\cal X$.
In fact, it is sufficient to construct the morphism 
on small neighbourhoods $U$ (in the sense of classical topology)
of each points.
By using a universal family on $U \times_T \cal X$, we get the morphism.
Since $s: \cal M \to T$ and $\Pic^0_{\cal X/T}\times_T \cal X \to T$ are
smooth over $T$,
it is sufficient to prove that
\begin{equation}
\frak a_{E_0}^*:H^1(\widehat{X} \times X,\Bbb Z)
 \to H^1(M_H(r,c_1,\Delta),\Bbb Z)
\end{equation}
is an isomorphism, if $X$ is a product of elliptic curves.
In order to prove this assertion, we shall show that 
\begin{equation}\label{eq:E0}
\frak a_{E_0}^*:\Pic^0(\widehat{X} \times X) \to \Pic^0(M(r,c_1,\Delta))
\end{equation}
is an isomorphism.
 Let $\cal E$ be a universal family on 
$M(r,c_1,\Delta)$.
For simplicity, we set $M:=M(r,c_1,\Delta)$.
Let $\widehat{X} \times X \to \Pic^0(X \times \widehat{X})$
be the isomorphism sending $(\hat{x},x) \in \widehat{X} \times X$
to $\cal P_{|\{\hat{x}\} \times X} \otimes \cal P_{|\widehat{X} \times \{x\}}$.
We set $\cal R:=\det p_{\widehat{X} \times M!}((\cal E-E_0 \otimes \cal O_M)
\otimes (\cal P-\cal O_{\widehat{X} \times X}))$.
By the construction of $\alpha_{E_0}$,
we get that $\cal R \cong 
(1_{\widehat{X}} \times \alpha_{E_0})^* \cal P \otimes L$,
where $L$ is the pull-back of a line bundle on $M$.
Since $\cal R_{|\{0\} \times M}\cong \cal O_M$, 
we get that $L\cong \cal O_{\widehat{X} \times M}$.
Hence we see that
\begin{equation}\label{eq:6-12-1}
\begin{split}
\alpha^*_{E_0}(\cal P_{|\{\hat{x}\}\times X}) &=
 \det p_{M!}((\cal E-E_0 \otimes \cal O_M)
\otimes(\cal P_{|\{\hat{x}\}\times X}-\cal O_X))\\
&=\det p_{M!}(\cal E\otimes(\cal P_{|\{\hat{x}\}\times X}-\cal O_X)).
\end{split}
\end{equation} 
In the same way, we see that
\begin{equation}\label{eq:6-12-2}
\begin{split}
\det\nolimits_{E_0}^*(\cal P_{|\widehat{X} \times \{x\}})&=
(\det \cal E \otimes \det E_0^{\vee} \otimes 
\det \cal E_{|M \times \{0\}}^{\vee})_{|M \times \{x\}}\\
&=\det p_{M!}(\cal E \otimes (k_x-k_0)),
\end{split}
\end{equation}
where $0 \in X$ is the zero of the group low.
In order to prove \eqref{eq:E0}, we shall consider the pull-backs of
$\alpha^*_{E_0}(\cal P_{|\{\hat{x}\}\times X})$
 and $\det_{E_0}^*(\cal P_{|\widehat{X} \times \{x\}})$
to $X^0 \times Y_0$. 

We denote the zero of $C_1$ and $C_2$ by $0_1$ and $0_2$ respectively.
For a point $q_k$ of $C_k$, $k=1,2$,
we set $l_k:=q_k-0_k$.
We also denote the pull-back of $l_k$ to $X=C_1 \times C_2$ by $l_k$. 
In the same way as in , we denote $\varpi_i^{!} (G)$, $i=0,1,\dots,n,a$ by
$G^i$, $G \in K(X)$.
We also denote $\cal O_X(D)^i$ by $\cal O_{X^0 \times Y_0}(D^i)$.
By simple calculations, we see that
\begin{equation}\label{eq:6-12-3}
\left\{
\begin{aligned}
& \det p_{X^0 \times Y_0!}(\cal E\otimes(\cal O_X(l_1)-\cal O_X)^a)
=\cal O_{X^0 \times Y_0}(dl^0_1-d_2\sum_{i=1}^n l_1^i)\\
& \det p_{X^0 \times Y_0!}(\cal E\otimes(\cal O_X(l_2)-\cal O_X)^a)
= \cal O_{X^0 \times Y_0}(\sum_{i=1}^n  l_2^i)\\
& \det p_{X^0 \times Y_0!}(\cal E\otimes(k_{(q_1,0_2)}-k_{(0_1,0_2)})^a)
= \cal O_{X^0 \times Y_0}(rl_1^0-r_2\sum_{i=1}^n l_1^i)\\
& \det p_{X^0 \times Y_0!}(\cal E\otimes(k_{(0_1,q_2)}-k_{(0_1,0_2)})^a)
=\cal O_{X^0 \times Y_0}(l_2^0).
\end{aligned}
\right.
\end{equation}
Since $d_2 r-d r_2=1$ and $\Pic^0(X^0 \times \Hilb_X^n) \cong
\Pic^0(X^0 \times Y_0)^{\frak S_n}$,
\eqref{eq:6-12-1}, \eqref{eq:6-12-2} and \eqref{eq:6-12-3} implies that
\eqref{eq:E0} holds.

We set 
\begin{equation}
K(r,c_1,\Delta):=\{\alpha \in K(X)|
\chi(\alpha \otimes E_0)=0, E_0 \in M_H(r,c_1,\Delta)\}.
\end{equation}
Let $\{U_i \}$ be an open covering of $M_H(r,c_1,\Delta)$
such that there are universal family $\cal F_i$ on each $U_i \times X$
and $\cal F_{i|(U_i \cap U_j) \times X} \cong 
\cal F_{j|(U_i \cap U_j) \times X}$.
Since the action of $\cal O_{U_i}^{\times}$ to
$\det p_{U_i!}(\cal F_i \otimes \alpha)$ is trivial,
we get a line bundle $\tilde{\kappa}(\alpha)$ on $M_H(r,c_1,\Delta)$.
Thus we obtain a homomorphism
\begin{equation}
\tilde{\kappa}:K(r,c_1,\Delta) \to \Pic(M_H(r,c_1,\Delta)).
\end{equation}
We note that there is a commutative diagram:
\begin{equation}
\begin{CD}
K(r,c_1,\Delta) @>{\tilde{\kappa}}>> \Pic(M_H(r,c_1,\Delta))\\
@V{\chr}VV @VV{c_1}V\\
H(r,c_1,\Delta) @>{\kappa_2}>> H^2(M_H(r,c_1,\Delta),\Bbb Z)
\end{CD}
\end{equation}
Let $K^2$ be the subgroup of $K(r,c_1,\Delta)$ generated by
$k_P-k_0$, $P \in X$
and $N$ the kernel of the Albanese map $K^2 \to X$.
Since $\ker (\chr)$ is generated by $\cal O_X(D)-\cal O_X$, 
$\cal O_X(D) \in \Pic^0(X)$ and $k_P-k_0$, $P \in X$,
$\eqref{eq:6-12-1}$ and $\eqref{eq:6-12-2}$ implies that
$\tilde{\kappa}$ induces an isomorphism $\ker(\chr)/N \to 
 \Pic^0(M_H(r,c_1,\Delta))$.
By using Theorem \ref{thm:H2} (3), we get the following theorem,
which is similar to [Y2, Thm. 0.1].
\begin{thm}\label{thm:Pic}
Under the same assumption as in Theorem \ref{thm:H2},
the following holds.\newline
$(1)$ $\frak a_{E_0}:M_H(r,c_1,\Delta) \to \widehat{X} \times X$
is an Albanese map.\newline
$(2)$
$\tilde{\kappa}:K(r,c_1,\Delta)/N \to \Pic(M_H(r,c_1,\Delta))$
is injective.\newline
$(3)$
$\Pic(M_H(r,c_1,\Delta))/\frak a_{E_0}^*(Pic(\widehat{X} \times X))$
is generated by $\tilde{\kappa}(K(r,c_1,\Delta))$.\newline
$(4)$ $\frak a_{E_0}^*(\Pic(\widehat{X} \times X)) 
\cap \tilde{\kappa}(K(r,c_1,\Delta)) \cong X \times \widehat{X}$.
\end{thm}

\section{appendix}
In this appendix, we shall show the following.

\begin{prop}\label{prop:1-24}
Let $L$ be an ample line bundle on $X$.
We assume that $\chi(L)=(c_1(L)^2)/2$ and $r$ are relatively prime.
Then $M_H(r,c_1(L),\Delta) \cong M_H(r,L,\Delta) \times \widehat{X}$,
where $M_H(r,L,\Delta)$ is the moduli space of determinant $L$.
In particular, $P(M_H(2,L,\Delta),z)=P(\Hilb_X^{2\Delta},z)$
for a general polarization $H$.
\end{prop}

\begin{pf}
For a stable sheaf $E \in M_H(r,c_1(L),\Delta)$,
$\lambda(E)$ denotes the point of $\widehat{X}$
which correspond to the line bundle $\det(E) \otimes L^{-1}$.
Let $\phi_L:X \to \widehat{X}$ be the morphism sending $x \in X$ to 
$T^*_x L \otimes L^{-1}$,
and $\varphi: \widehat{X} \to X$ the morphism such that 
$\phi_L \circ \varphi=n^2_{\widehat{X}}$,
where $T_x: X \to X$ is the translation defined by $x$
and $n^2=\chi(L)^2=\deg \phi_L$.
Since $(r,n^2)=1$, there are integers $k$ and $k'$ such that $rk+n^2 k'=1$.
We denote the Poincar\'{e} line bundle on $X \times \widehat{X}$ 
by $\cal P$.
Let $A:M_H(r,c_1(L),\Delta) \to M_H(r,L,\Delta)\times \widehat{X}$
be the morphism sending $F \in M_H(r,c_1(L),\Delta)$
to $(T^*_{-k'\varphi \circ \lambda(F)}
(F \otimes \cal P_{-k\lambda(F)}),
\lambda(F))$
and $B:M_H(r,L,\Delta)\times \widehat{X} \to M_H(r,c_1(L),\Delta)$
the morphism sending
$(E,x) \in M_H(r,L,\Delta)\times \widehat{X}$ to 
$T^*_{k'\varphi(X)}E \otimes \cal P_{kx}$.
For an element $(E,x)$ of $M_H(r,L,\Delta)\times \widehat{X}$,
$\det(T^*_{k'\varphi(x)}E \otimes \cal P_{kx})
\cong T_{k'\varphi(x)}^*L \otimes \cal P_{rkx}
\cong L \otimes \cal P_{k'\phi_L \circ \varphi (x)} \otimes \cal P_{rkx}
\cong L \otimes \cal P_{(n^2k'+rk)x}=L \otimes \cal P_{x}$.
Hence $\lambda \circ B((E,x))=x$.
Then it is easy to see that $A \circ B$ and $B \circ A$ are identity 
morphisms.
Hence $A:M_H(r,c_1(L),\Delta) \to M_H(r,L,\Delta) \times \widehat{X}$
is an isomorphism.
\end{pf}

Let $\Bbb D(X)$ and $\Bbb D(\widehat{X})$ be the derived categories of
$X$ and $\widehat{X}$ respectively.
Let $\cal S:\Bbb D(X) \to \Bbb D(\widehat{X})$ be the Fourier-Mukai
transform [Mu4].
Then the morphism $\alpha:=\alpha_{E_0}$ defined in {\bf 3.4} satisfies that
$\alpha(E)=\det \cal S(E) \otimes (\det\cal S(E_0))^{-1}$.
Thus $\alpha_{E_0}$ is also defined by Fourier-Mukai transform.
By using [M4], we shall treat the case $2r\Delta=2$
(at least, Mukai treated the case where 
$X$ is a principally polarized Abelian surface). 
\begin{prop}
Let $L$ be an ample divisor.
If $2r\Delta=2$, then for a general polarization $H$,
$\alpha:M_H(r,L,\Delta) \to X$ is an isomorphism.
\end{prop}

\begin{pf}
Since $rc_2-(r-1)(L^2)/2=1$ and $\chi(L)=(L^2)/2$,
$r$ and $\chi(L)$ are relatively prime.
We shall choose an element $E$ of $M_H(r,L,\Delta)$ and let
$\xi:X \times \widehat{X} \to M(r,c_1(L),\Delta)$ be the morphism
sending $(x,y) \in X \times \widehat{X}$ to $T_x^*E \otimes \cal P_y$.
Then $\lambda \circ \xi(x,y)=\phi_L(x)+ry$.
Let $f:X \to X \times \widehat{X}$ be the morphism such that 
$f(x)=(rx,-\phi_L(x))$.
Since $\# \ker \phi_L=\chi(L)^2$ and $r$ are relatively prime,
$f$ is injective.
Let $g:\widehat{X} \to X \times \widehat{X}$ be the morphism
such that $g(y)=(k'\varphi(y),ky)$.
Then $f \times g:X \times \widehat{X} \to X \times \widehat{X}$
is an isomorphism.
In fact, if $(rx+k'\varphi(y), -\phi_L(x)+ky)=(0,0)$,
then $\phi_L(rx+k'\varphi(y))=r\phi_L(x)+n^2k'y=0$.
Hence $y=(n^2k'+rk)y=0$. Since $f$ is injective,
$x=0$,
which implies that $f \times g$ is injective.
Therefore $f \times g$ is an isomorphism.
Then we get a morphism
$\xi \circ f: X \to M(r,L,\Delta)$.
Replacing $E$ by $E \otimes L^{\otimes m}$,
we may assume that there is an exact sequence
$0 \to \cal O_X^{\oplus (r-1)} \to E \to I_Z \otimes L \to 0$,
where $I_Z$ is the ideal sheaf of a codimension 2
subscheme $Z$ of $X$.
By our assumption on Chern classes,
$1/r=\Delta(E)=\deg Z-(r-1)/r\chi(L)$.
For simplicity,
we denote $\det \cal S (?)$ by $\delta(?)$.
Then we see that
$\delta(T^*_x E \otimes \cal P_y)
=\delta(I_{T_{-x}(Z)} \otimes T^*_xL \otimes \cal P_y)
=\delta(I_{Z-(\deg Z)x} \otimes L \otimes \cal P_{\phi_L(x)+y})
=\delta(L \otimes \cal P_{\phi_L(x)+y})\otimes \cal P_{-Z+(\deg Z)x}
=\det T^*_{\phi_L(x)+y}(\cal S(L)) \otimes 
\cal P_{-Z+(\deg Z)x} 
=\delta (L) \otimes \cal P_{\phi_{\delta(L)}(\phi_L(x)+y)+(\deg Z)x-Z}$.
Hence $\alpha \circ \xi \circ f(x)
=\alpha \circ \xi \circ f(0)+(r-1)\phi_{\delta(L)}\circ \phi_L(x)+r(\deg Z)x$.
By the proof of [Mu4, Prop.1.23],
$\phi_{\delta(L)}(\phi_L(x))=-\chi(L)x$.
Since $r\deg Z=1+(r-1)\chi(L)$,
we get that $\alpha\circ \xi \circ f(x)=\alpha \circ \xi \circ f(0)+x$.
Thus $\alpha\circ \xi \circ f(x)$ is an isomorphism.
Therefore we get that $\alpha:M_H(r,L,\Delta) \to X$ is an isomorphism.
\end{pf}

\end{document}